\documentclass[a4paper, 11pt, twoside]{article}
\usepackage{mfpaperstuff}
\usepackage[normalem]{ulem}
\usepackage{marvosym}

\begin{document}
\newcommand\lir[3]{\operatorname c^{#1}_{#2#3}}
\newcommand\sch[1]{\cals(#1,#1)}
\renewcommand{\crefrangeconjunction}{--}
\medmuskip=3mu plus 1mu minus 1mu
\newcommand\hb[1]{\hbox to .5cm{\hfil#1\hfil}}
\newcounter{casy}
\newcounter{caseone}
\newcommand\imo{i-1}
\newcommand\e{\mathrm{e}}
\newcommand\f{\mathrm{f}}
\newcommand\po{+\negthinspace1}
\newcommand\bsm{\begin{smallmatrix}}
\newcommand\esm{\end{smallmatrix}}
\newcommand\pon[2]{#1\circ#2}
\newcommand\ptw[2]{#1\ast#2}
\newcommand{\aaa}{\mathfrak{A}_}
\newcommand\lad[2]{\operatorname{rmp}_{#1}(#2)}
\newcommand\alad[2]{\operatorname{rmp}^+_{#1}(#2)}
\newcommand\rlad[2]{\operatorname{rmp}^-_{#1}(#2)}
\newcommand\reall[2]{#2^{\!\triangledown#1}}
\newcommand\addall[2]{#2^{\!\triangle#1}}
\newcommand\renor[2]{#2^{\!\blacktriangledown#1}}
\newcommand\rest{^{\operatorname{rest}}}
\newcommand\zpz{\bbz/p\bbz}
\newcommand\nor[2]{\operatorname{nor}_{#1}(#2)}
\newcommand\remo[2]{\operatorname{rem}_{#1}(#2)}
\newcommand\spe[1]{\ensuremath{\operatorname{S}^{#1}}}
\newcommand\sid[1]{\operatorname{D}_{#1}}
\newcommand\weyl[1]{\ensuremath{\Delta^{#1}}}
\newcommand\schursimp[1]{\operatorname{L}^{#1}}
\newcommand\yper[1]{\operatorname{M}^{#1}}
\newcommand\tht[1]{\Theta_{#1}}
\newcommand\hth[1]{\hat\Theta_{#1}}
\newcommand\ho[2]{\operatorname{Hom}_{\bbf\sss n}(#1,#2)}
\newcommand\noregdown[2]{\lightning_{#1}(#2)}
\newcommand\soc[1]{\operatorname{soc}(#1)}
\newcommand\mullna[1]{\operatorname{m}_{#1}}
\newcommand\mull[2]{\mullna{#1}(#2)}
\newcommand\mul[1]{\mull3{#1}}
\newcommand\sgn{\operatorname{sgn}}
\newcommand\mood{\textup{-mod}}
\newcommand\res{\operatorname{res}}
\newcommand\emax[1]{\e_{#1}^{{(\max)}}}
\newcommand\ep{\epsilon}
\newcommand\eps[2]{\ep_{#1}#2}
\newcommand\dn[2]{[\spe{#1}\,{:}\,\sid{#2}]}
\newcommand\len[1]{#1'_1}
\newcommand\roof[1]{\operatorname{roof}(#1)}
\newcommand\cp[3]{\operatorname{CP}^{#1}_{#2}(#3)}
\newcommand\rtab[1]{\operatorname{Re}(#1)}
\newcommand\stre{^{\circ}}
\newcommand\ins[2]{{}^{\Rsh}_{#1}#2}
\newcommand\insm[1]{#1^+}
\newcommand\muti[2]{\{#1\}^{#2}}
\newcommand\opm{\muti1{p-1}}
\newcommand\row{\stackrel{\operatorname{row}}\longleftrightarrow}
\newcommand\col{\stackrel{\operatorname{col}}\longleftrightarrow}
\newcommand\csgn[2]{\operatorname{colsign}_{#1#2}}
\newcommand\rec[3]{#1[#2,#3]}
\newcommand\full[1]{\operatorname{full}(#1)}
\newcommand\wei[2]{\operatorname{wt}(#2)}
\newcommand\twofac{\ensuremath{\cala}}
\newcommand\jmset{\ensuremath{\calj}}
\newcommand\al[1]{\ensuremath{\reall{\if 0#1\else\pm\fi#1}\la}}
\newcommand\pt h
\newcommand\quo[2]{#2^{(#1)}}
\newcommand\oquo[2]{#2^{[#1]}}
\newcommand\opq{ordered $p$-quotient}
\newcommand\opqd[1]{\big[\oquo0{#1},\dots,\oquo{p-1}{#1}\big]}
\newcommand\pqd[1]{\big(\quo0{#1},\dots,\quo{p-1}{#1}\big)}
\newcommand\nrd[1]{\ensuremath{\noregdown{#1}\la}}
\newcommand\rone{R-partition of type~I}
\newcommand\rones{R-partitions of type~I}
\newcommand\rtwo{R-partition of type~II}
\newcommand\rtwos{R-partitions of type~II}
\newcommand\cf{composition factor}
\newcommand\qs{$p$-quotient-separated}
\renewcommand\sc{self-conjugate}
\newcommand\scy{self-conjugacy}
\newcommand\jm{JM-partition}
\title{The irreducible representations of the alternating\\group which remain irreducible in characteristic $p$}
\runninghead{Representations of the alternating group which remain irreducible in characteristic $p$}
\toptitle

\begin{abstract}
Let $p$ be an odd prime, and $\aaa n$ the alternating group of degree $n$. We determine which ordinary irreducible representations of $\aaa n$ remain irreducible in characteristic $p$, verifying the author's conjecture from \cite{altred}. Given the preparatory work done in [\textit{op.\ cit.}], our task is to determine which self-conjugate partitions label Specht modules for the symmetric group in characteristic $p$ having exactly two composition factors. This is accomplished through the use of the Robinson--Brundan--Kleshchev `$i$-restriction' functors, together with known results on decomposition numbers for the symmetric group and additional results on the Mullineux map and homomorphisms between Specht modules.
\end{abstract}


\renewcommand\baselinestretch{1.006}
\Yboxdim{12pt}

\section{Introduction}

An interesting question for any finite group is to ask which ordinary irreducible representations of that group remain irreducible in characteristic $p$. For the symmetric group $\sss n$ this amounts to classifying the irreducible Specht modules, and this problem was solved several years ago, through the combined efforts of James, Mathas, Lyle and the author \cite{jm1,jmp2,slred,mfred,mfirred}; when $p$ is odd, the irreducible Specht modules are precisely those labelled by \emph{\jm s}. In this paper we address the case of the alternating group $\aaa n$. The author considered this problem in \cite{altred}, solving it completely in the case $p=2$ and presenting a conjectured solution for odd $p$, which we prove here.

As with many problems concerning the representation theory of the alternating group, our technique is to translate the problem to one about the symmetric group, using elementary Clifford theory to transfer results between the two settings. For the problem at hand, this translation was done in \cite{altred}, where the main problem for the alternating groups was reduced to the question of which Specht modules labelled by self-conjugate partitions have exactly two composition factors (with multiplicity). So the present paper is concerned entirely with the representation theory of the symmetric group.

Our main result when $p\gs5$ is that the self-conjugate partitions labelling Specht modules with composition length $2$ are the partitions which we called \emph{R-partitions} in \cite{altred}; these have a simple description in terms of hook lengths in the Young diagram. The fact that the corresponding Specht modules have composition length $2$ was shown in \cite{altred}, so the task undertaken in this paper is to prove the converse. The same applies in the case $p=3$, where the classification takes a slightly different form.

Our basic strategy involves applying the results of Brundan and Kleshchev concerning Robinson's $i$-restriction functors $\e_i$. We suppose $\la$ is a self-conjugate partition which is not an R-partition or a \jm; then we must show that the Specht module $\spe\la$ has at least three composition factors. By removing all the removable nodes of residue $0$ from the Young diagram of $\la$, one obtains a self-conjugate partition $\reall0\la$, and the Brundan--Kleshchev results imply that the composition length of $\spe\la$ is at least that of $\spe{\reall0\la}$; so by induction, we may assume that either $\reall0\la$ is an R-partition or a \jm, or $\reall0\la=\la$. Similarly, for any $i\neq0$ we can define a self-conjugate partition $\al i$ by repeatedly removing removable nodes of residues $i$ and $p-i$ from $\la$, and we can make the same inductive assumption about $\al i$. This restricts the possibilities for $\la$ considerably. In fact, we can strengthen this inductive argument using James's Regularisation Theorem, which gives an explicit composition factor $\sid{\la\rest}$ of the Specht module $\spe\la$, for any $\la$. It is very helpful for our purposes to be able to tell, for a given partition $\la$, when there is an $r$ such that $\e_i^r\spe\la\neq0$ but $\e_i^r\sid{\la\rest}=0$. An important new result in this paper (\cref{renorl}) is an explicit combinatorial criterion for this; it is surprising to the author that this does not seem to have been discovered before.

This inductive argument deals with most cases. Several of the remaining cases are eliminated with using the theory of \emph{Rouquier blocks}, whose decomposition numbers are very well understood. By developing the associated abacus combinatorics we exploit these results, together with the classification of irreducible Weyl modules, to show that our main theorem holds for the so-called \emph{\qs} partitions.

These arguments leave just one family of Specht modules to deal with, which we approach using the theory of homomorphisms between Specht modules. Establishing the existence of a non-zero homomorphism $\spe\la\to\spe\nu$ shows that $\spe\la$ and $\soc{\spe\nu}$ share a composition factor, and we use this in a certain special case to show that $\spe\la$ has at least three composition factors. The Specht homomorphism that we use is constructed as the composition of two well-known homomorphisms, namely the \emph{one-node Carter--Payne homomorphism} and the \emph{regularisation homomorphism}. However, we have considerable work to do in showing that the composition is non-zero in our particular situation. To do this, we give a result describing the least dominant tableau occurring when the regularisation homomorphism is expressed in terms of semistandard homomorphisms; as a by-product, this gives a new proof of the (non-trivial) fact that the regularisation homomorphism is non-zero.

This homomorphism result, together with a small lemma concerning the \emph{Mullineux map} (which describes the effect of the functor $-\otimes\sgn$ on simple modules), is enough to complete the proof. We conclude the paper with a simple corollary which shows that the only irreducible representations of $\aaa n$ remaining irreducible modulo every prime are the one-dimensional representations.

\begin{acks}
The author was inspired to re-visit this problem at an Oberwolfach mini-workshop on the representation theory of the symmetric groups in 2011; he is very grateful to Suzanne Danz and David Hemmer for the invitation to this workshop, and to David Hemmer for subsequent discussions on this problem.

The author is also indebted to the referee for a very careful reading of the paper.
\end{acks}

\section{Representation theory of the symmetric group}\label{backsec}

In this section, we recall some essential background on the representation theory of the symmetric group. Throughout this paper $n$ is a non-negative integer and $\bbf$ is a field of characteristic $p$; we use the convention that the characteristic of a field is the order of its prime subfield, so $p\in\{2,3,5,\dots\}\cup\{\infty\}$.

\subsection{Partitions and Specht modules}\label{partnsubsec}

A \emph{composition} of $n$ is a sequence $\la=(\la_1,\la_2,\dots)$ of non-negative integers summing to $n$. If $n$ is not specified, we write $|\la|$ for the sum of the terms of $\la$. When writing compositions, we usually omit trailing zeroes and group together consecutive equal parts with a superscript; the unique composition of $0$ is denoted $\varnothing$. A composition which is weakly decreasing is called a \emph{partition}. We often identify a composition $\la$ with its \emph{Young diagram}, which is the set
\[
\rset{(r,c)\in\bbn^2}{c\ls\la_r}.
\]
Elements of the Young diagram of $\la$ are called \emph{nodes} of $\la$; more generally, a \emph{node} is any element of $\bbn^2$. We adopt the English convention for drawing Young diagrams, in which $\la$ is drawn with left-justified rows of boxes of lengths $\la_1,\la_2,\dots$ successively down the page.

If $\la$ is a partition, the \emph{conjugate} partition $\la'$ is defined by
\[
\la'_i=\left|\rset{j\in\bbn}{i\ls\la_j}\right|,
\]
or, in terms of Young diagrams, by reflecting along the main diagonal. $\la$ is \emph{self-conjugate} if $\la'=\la$. $\la$ is \emph{$p$-restricted} if $\la_i-\la_{i+1}<p$ for all $i\gs1$, and \emph{$p$-regular} if $\la'$ is $p$-restricted.

A node $(r,c)$ of a partition $\la$ is \emph{removable} if it can be removed from $\la$ to leave a smaller partition (i.e.\ if $c=\la_r>\la_{r+1}$), while a node $(r,c)$ not in $\la$ is an \emph{addable node} of $\la$ if it can be added to $\la$ to give a larger partition.

The \emph{$p$-residue} of a node $(r,c)$ is the residue of $c-r$ modulo $p$ (or simply the integer $c-r$, when $p=\infty$). If a node has residue $i$, we call it an \emph{$i$-node}.

If $\la$ is a composition of $n$, let $\yper\la$ denote the \emph{Young permutation module} for $\bbf\sss n$ corresponding to $\la$, as defined in \cite[\S4]{jbook}. If $\la$ is a partition, let $\spe\la$ denote the \emph{Specht module} corresponding to $\la$. If $\la$ is $p$-restricted, then $\spe\la$ has a simple socle $\sid\la$, and the modules $\sid\la$ afford all the irreducible representations of $\bbf\sss n$ as $\la$ ranges over the set of $p$-restricted partitions of $n$. When $p=\infty$, we have $\sid\la=\spe\la$, so the characters $\chi^\la$ of the Specht modules give all the ordinary irreducible characters of $\sss n$.

\begin{rmk}
It is slightly more traditional to label the simple $\bbf\sss n$-modules by $p$-regular partitions: if $\la$ is $p$-regular, then $\spe\la$ has a simple cosocle $\operatorname{D}^\la$, and these modules also afford all the irreducible representations of $\bbf\sss n$. It is well known how to convert from one convention to the other; we have chosen the $p$-restricted convention in this paper because it aligns better with some of the references that we cite.
\end{rmk}

We shall also briefly need to consider Weyl modules for the Schur algebra $\sch n$, for which we refer to the book by Green \cite{gbook}. We let $\weyl \la$ denote the \emph{Weyl module} (also called the \emph{standard module}) labelled by the partition $\la$ of $n$, and $\schursimp\la$ its unique irreducible quotient.

\subsection{James's Regularisation Theorem}

The main aim in this paper is to consider Specht modules with very few composition factors. A very helpful fact in this endeavour when $p<\infty$ is that we know an explicit composition factor of every Specht module. This result is James's Regularisation Theorem, which we phrase here in terms of $p$-restricted partitions.

Suppose $l\gs0$. Define the $l$th \emph{ramp} in $\bbn^2$ to be the set of nodes $(r,c)$ for which $c-1+(p-1)(r-1)=l$. If $l<m$, we say that ramp $m$ is \emph{later} than ramp $l$. If $\la$ is a partition, the \emph{$p$-restrictisation} of $\la$ is the $p$-restricted partition $\la\rest$ obtained by moving all the nodes in each ramp as far to the left within that ramp as possible. ($\la\rest$ is simply called the \emph{$p$-restriction} of $\la$ in \cite{flm}, but we introduce the slightly absurd term \emph{restrictisation} here to avoid confusion with restriction in the sense of restricting to subgroups, which we shall consider a great deal. The linguistically sensitive reader may rest assured that we shall use this term as little as possible.)

\begin{eg}
Take $\la=(8,6,2,1^2)$ and $p=3$. Then $\la\rest=(6,5,4,2,1)$, as we can see from the following Young diagrams, in which we label each node with the number of the ramp in which it lies.
\[
\young(01234567,234567,45,6,8),\qquad\young(012345,23456,4567,67,8).
\]
\end{eg}

Now we can state part of James's Regularisation Theorem, translated to the $p$-restricted convention.

\begin{thmc}{j1}{Theorem A}\label{jrest}
Suppose $\la$ is a partition. Then $\dn\la{\la\rest}=1$.
\end{thmc}

\subsection{The $p$-core and $p$-weight of a partition}

If $\la$ is a partition, the \emph{rim} of $\la$ is defined to be the set of all nodes $(r,c)$ of $\la$ such that $(r+1,c+1)$ is not a node of $\la$. Given any node $(r,c)$ of $\la$, the \emph{$(r,c)$-rim hook} of $\la$ is the connected portion of the rim running from the node $(r,\la_r)$ down to the node $(\la'_c,c)$. The \emph{$(r,c)$-hook length} of $\la$ is the number of nodes in the $(r,c)$-rim hook, i.e.\ $\la_r-r+\la'_c-c+1$. If $p<\infty$, then we say that $\la$ is a \emph{$p$-core} if none of the hook lengths of $\la$ is divisible by $p$, or equivalently if none of the hook lengths equals $p$.

If $\la$ is an arbitrary partition, the \emph{$p$-core of $\la$} is obtained as follows. Choose a node $(r,c)$ of $\la$ such that the $(r,c)$-hook length equals $p$, and delete the $(r,c)$-rim hook from $\la$; repeat until a $p$-core is obtained. This $p$-core is independent of the choice of rim hook deleted at each stage, and hence so is the number of rim hooks deleted; this number is called the \emph{$p$-weight} of $\la$. 

\subsection{The sign representation}\label{signsubsec}

Let $\sgn$ denote the one-dimensional sign representation of $\sss n$. This gives rise to a functor $-\otimes\sgn:\bbf\sss n\mood\to\bbf\sss n\mood$, which takes simple modules to simple modules. The effect of this functor on Specht modules is well-known; let $M^\ast$ denote the dual of a module $M$.

\begin{thmc}{jbook}{Theorem 8.15}\label{815}
If $\la$ is a partition, then $\spe\la\otimes\sgn\cong(\spe{\la'})^\ast$.
\end{thmc}

Since every simple $\bbf\sss n$-module is self-dual \cite[Theorem 11.5]{jbook} and every Specht module is indecomposable when $p\gs3$ \cite[Corollary 13.18]{jbook}, \cref{jrest,815} imply the following.

\begin{lemmac}{altred}{Proof of Theorem 4.1}\label{scmull}
Suppose $p\gs3$, and $\la$ is a self-conjugate partition such that $\spe\la$ has exactly two composition factors $S,T$. Then $S\otimes\sgn\cong T\ncong S$.
\end{lemmac}

Now we consider simple modules. Since $S\otimes\sgn$ is simple whenever $S$ is, there is an involution $\mullna p$ on the set of $p$-restricted partitions of $n$, defined by $\sid\la\otimes\sgn\cong\sid{\mull p\la}$. This involution is known as the \emph{Mullineux map}, since Mullineux \cite{mull} gave (albeit without proof) the first of several known recursive combinatorial descriptions of the map. We do not give this algorithm here, since for this paper we just need the following simple result concerning $\mullna p$.

\begin{lemma}\label{mulllemma}
Suppose $\la$ is a $p$-restricted partition. Then $\len\la+\len{\mull p\la}\nequiv1\pmod p$.
\end{lemma}

\begin{pf}
We use the work of Ariki et al.\ \cite{akt} which addresses the relationship between the labellings of simple modules by $p$-restricted partitions and by Littelmann paths. Following Kreiman et al.\ \cite{klmw}, they define for each $p$-restricted partition $\la$ a $p$-core $\roof\la$, which has the following two properties:
\begin{description}
\item[{\rm\cite[Lemma 2.4(3)]{akt}}] $\len{\roof\la}=\len\la$;
\item[{\rm\cite[Proposition 5.21]{akt}}] $\roof{\mull p\la}=\roof\la'$.
\end{description}
From these properties we have $\len\la+\len{\mull p\la}=\len{\roof\la}+\roof\la_1$. This cannot be congruent to $1$ modulo $p$, since then the $(1,1)$-hook length of $\roof\la$ would be divisible by $p$, contradicting the fact that $\roof\la$ is a $p$-core.
\end{pf}

It would be very interesting to see a more direct proof of \cref{mulllemma} using Mullineux's algorithm, for example.

\subsection{\jm s}\label{jmsec}

Now we describe the partitions which label irreducible Specht and Weyl modules when $p\gs3$. Define the \emph{$p$-power diagram} of a partition $\la$ to be the diagram obtained by filling the $(r,c)$-box in the Young diagram of $\la$ with the $p$-adic valuation of the $(r,c)$-hook length, for each node $(r,c)$ of $\la$. Say that $\la$ is a \emph{$p$-\jm} (or simply a \emph{\jm}) if the following property holds: every non-zero entry in the $p$-power diagram is either equal to all the other entries in the same row or equal to all the other entries in the same column.

\begin{eg}
The partition $(19,11,2^3,1^2)$ is a $3$-\jm, as we see from its $3$-power diagram:
\[
\young(0020010010000100100,00200100100,11,00,00,0,0).
\]
\end{eg}

\begin{thmciting}{\textup{\textbf{\cite{jm1,slred,mfred,mfirred}.}}\ }\label{irredspecht}
If $p\gs3$, then the Specht module $\spe\la$ is irreducible if and only if $\la$ is a $p$-JM partition.
\end{thmciting}

A similar result applies for irreducible Weyl modules (without the restriction $p\gs3$). Observe that a \jm{} is $p$-restricted if and only if every non-zero entry in its $p$-power diagram is equal to every entry in the same row. A similar statement applies to $p$-regular \jm s (which are also known as \emph{Carter partitions}).

\begin{thmc}{jm1}{Theorem 4.5}
The Weyl module $\weyl\la$ is irreducible if and only if $\la$ is a $p$-restricted \jm.
\end{thmc}

\subsection{R-partitions}

We now describe a family of self-conjugate partitions introduced in \cite{altred} which label irreducible representations of the alternating group that remain irreducible in characteristic $p$.

Given a partition $\la$, construct the $p$-power diagram of $\la$ as above. Say that $\la$ is an \emph{R-partition} if $\la$ is \sc{} and there is a \emph{distinguished node} $(r,r)$ of $\la$ such that:
\begin{itemize}
\item
the $(r,r)$-entry in the $p$-power diagram is non-zero, and
\item
any non-zero entry in the $p$-power diagram \emph{other than the $(r,r)$-entry} is either equal to all the entries in its row or equal to all the entries in its column.
\end{itemize}

R-partitions were studied in detail in \cite{altred}; later, we shall cite some results describing abacus displays for R-partitions. For now, we recall the two distinct types of R-partitions described in \cite[\S4.2]{altred}. Suppose $\la$ is an R-partition, with distinguished node $(r,r)$.

\begin{enumerate}
\item
$\la$ is an \emph{\rone} if $r=1$. In this case, removing the $(1,1)$-rim hook from $\la$ leaves a self-conjugate $p$-core $\xi$ with $\xi_1\ls\frac12(p-1)$.
\item
$\la$ is a \emph{\rtwo} if the $(r,r)$-hook length of $\la$ equals $p$. In this case, removing the $(r,r)$-rim hook from $\la$ leaves a self-conjugate \jm.
\end{enumerate}
Note that the \rtwos{} include all self-conjugate partitions of $p$-weight $1$. It was shown in \cite{altred} that every R-partition is either of type I or type II; however, the two types are not mutually exclusive. 

\begin{eg}
Take $p=5$. Then the partitions $(13,3^2,1^{10})$ and $(14,10,5,4,3,2^5,1^4)$ are R-partitions of types I and II respectively, with distinguished nodes $(1,1)$ and $(3,3)$, as we see from their $5$-power diagrams.
\[
\young(2001000010000,000,000,1,0,0,0,0,1,0,0,0,0)\qquad\young(00000100000000,0000010000,00100,0000,000,11,00,00,00,00,0,0,0,0)
\]
\end{eg}

\section{The alternating group and the main result}\label{altsec}

In this section we introduce the alternating group and give our main result, which we then re-cast in the symmetric group setting.

Suppose throughout this section that $\bbf$ is a splitting field for $\aaa n$, and that $p\neq2$. We also assume that $n\gs2$ (so that $\aaa n$ has index $2$ in $\sss n$), the case $n=1$ being trivial. If $M$ is an $\bbf\sss n$-module and $H\ls\sss n$, let $\res_HM$ denote the restriction of $M$ to $H$. If $M$ is irreducible, then by basic Clifford theory $\res_{\aaa n}M$ is irreducible if $M\otimes\sgn\ncong M$, and otherwise $\res_{\aaa n}M$ splits as the direct sum $M^+\oplus M^-$ of two irreducible modules. Furthermore, all irreducible $\bbf\aaa n$-modules arise in this way.

If $p=\infty$ then, by \cref{815} and the fact that simple $\bbf\sss n$-modules are self-dual, we have $\spe\la\otimes\sgn\cong\spe{\la'}$, so the character $\chi^\la$ of $\spe\la$ restricts to an irreducible character $\psi^\la$ of $\aaa n$ if $\la\neq\la'$ (and in this case $\chi^{\la'}$ also restricts to $\psi^\la$), while if $\la=\la'$ then $\chi^\la$ restricts to the sum of two irreducible characters $\psi^{\la+},\psi^{\la-}$. With this notation, we can give our main result.

\bigskip
\begin{mdframed}[innerleftmargin=3pt,innerrightmargin=3pt,innertopmargin=3pt,innerbottommargin=3pt,roundcorner=5pt,innermargin=-3pt,outermargin=-3pt]
\begin{thm}\label{mainthm}
Suppose $\bbf$ is a splitting field for $\aaa n$ of characteristic $p$, and $\psi=\psi^\la$ or $\psi^{\la\pm}$ is an ordinary irreducible character of $\aaa n$.
\begin{enumerate}
\item
If $p=3$, then $\psi$ is irreducible over $\bbf$ if and only if one of the following holds:
\begin{enumerate}
\item
$\la$ is a $3$-\jm;
\item
$\la$ has $3$-weight $1$;
\item
$\la$ is an \rone;
\item
$\la=(3^3)$.
\end{enumerate}
\item
If $p\gs5$, then $\psi$ is irreducible over $\bbf$ if and only if $\la$ is a $p$-\jm{} or an R-partition.
\end{enumerate}
\end{thm}
\end{mdframed}
\medskip

Most of \cref{mainthm} has already been proved in \cite{altred}, beginning with the following reduction of the problem to the representation theory of $\sss n$.

\begin{thmc}{altred}{Theorem 4.1}\label{altred4.1}
Suppose $\bbf$ is a splitting field for $\aaa n$ of odd characteristic $p$, and $\psi=\psi^\la$ or $\psi^{\la\pm}$ is an irreducible character of the alternating group $\aaa n$. Then $\psi$ is irreducible over $\bbf$ if and only if one of the following holds.
\begin{enumerate}
\item
$\la$ is a $p$-JM partition.
\item
$\la=\la'$ and the Specht module $\spe\la$ has exactly two composition factors.
\end{enumerate}
\end{thmc}

So from now on we can restrict attention entirely to the symmetric group, and prove the following.

\medskip
\begin{mdframed}[innerleftmargin=3pt,innerrightmargin=3pt,innertopmargin=3pt,innerbottommargin=3pt,roundcorner=5pt,innermargin=-3pt,outermargin=-3pt]
\begin{thm}\label{main}
Suppose $\bbf$ has characteristic $p$, and $\la$ is a self-conjugate partition. Then $\spe\la$ has exactly two composition factors if and only if one of the following holds.
\begin{itemize}
\item
$\la$ has $p$-weight $1$.
\item
$\la$ is an \rone.
\item
$p\gs5$ and $\la$ is an \rtwo.
\item
$p=3$ and $\la=(3^3)$.
\end{itemize}
\end{thm}
\end{mdframed}
\bigskip

We have already proved the `if' part of \cref{main} in \cite[\S5.1]{altred}. So the remainder of this paper is dedicated to proving the `only if' part. We do this in \cref{pfsec}, after we have recalled some more background and developed further tools.

\section{Restriction functors}\label{resfunc}

In this section, we describe some results on restriction functors, which will be our main tool. The definition of these functors goes back to Robinson \cite{gdbr}, though our main reference here is the survey of Brundan and Kleshchev \cite{bk}. We translate the partition combinatorics from \cite{bk} to the $p$-restricted convention.

\subsection{The restriction functors $\e_i$}\label{eisubsec}

Assume throughout this section that $p<\infty$. We shall feel free to identify $\zpz$ with the set $\{0,\dots,p-1\}$. In \cite[\S2.2]{bk}, Brundan and Kleshchev introduce the \emph{$i$-restriction operators} $\e_i:\bbf\sss n\mood\to\bbf\sss{n-1}\mood$, for $i\in\zpz$. These are exact functors, and have the property that
\[
\res_{\sss{n-1}}M\cong\bigoplus_{i\in\zpz}\e_iM
\]
for any $\bbf\sss n$-module $M$. In fact, if $M$ lies in a single block of $\bbf\sss n$ then the non-zero $\e_iM$ are precisely the block components of $\res_{\sss{n-1}}M$. In addition, $\e_i$ and $\e_j$ commute unless $j=i\pm1$.

The functors $\e_i$ are defined for all $n>0$, so it makes sense to define powers $\e_i^r$ for $r\gs0$. In fact, it is possible to define \emph{divided powers} $\e_i^{(r)}$, with the property \cite[Lemma 2.6]{bk} that
\[
\e_i^rM\cong\bigoplus_{k=1}^{r!}\e_i^{(r)}M
\]
for any $M$. If $M$ is non-zero, we define $\eps iM=\max\rset{r\gs0}{\e_i^{(r)}M\neq0}$, and $\emax iM=\e_i^{(\eps iM)}M$.

Now we recall some results describing the effect of these operators on Specht modules and simple modules. Given a partition $\la$ and $i\in\zpz$, let $\remo i\la$ denote the number of removable $i$-nodes of $\la$ and $\reall i\la$ the partition obtained by removing all the removable $i$-nodes from $\la$. Then the following is just a refinement of the classical Branching Rule for Specht modules.

\begin{lemma}\label{spechtrest}
Suppose $\la$ is a partition, and $i\in\zpz$. Then $\eps i{\spe\la}=\remo i\la$, and $\emax i\spe\la\cong\spe{\reall i\la}$.
\end{lemma}

\begin{eg}
Suppose $p=3$ and $\la=(4,3,1^2)$. The residues of the removable nodes of $\la$ are as follows:
\Yinternals0
\[
\yngres(3,4,3,1^2).
\]
Hence we have
\begin{alignat*}2
\eps0\spe\la&=2,&\qquad \emax0\spe\la&\cong\spe{(3^2,1)},\\
\eps1\spe\la&=1,&\qquad \emax1\spe\la&\cong\spe{(4,2,1^2)},\\
\eps2\spe\la&=0,&\qquad \emax2\spe\la&=\spe\la.
\end{alignat*}
\end{eg}

The corresponding result for the simple modules $\sid\la$ is more complicated, and requires some terminology. Define a \emph{sign sequence} to be a finite string of $+$ and $-$ signs, and define the \emph{reduction} of a sign sequence to be the sequence obtained by successively deleting adjacent pairs $-+$. If $\la$ is a $p$-restricted partition, define the \emph{$i$-signature} of $\la$ to be the sign sequence obtained by examining the addable and removable $i$-nodes of $\la$ from top to bottom, writing a $+$ for each addable $i$-node and a $-$ for each removable $i$-node, and define the \emph{reduced $i$-signature} of $\la$ to be the reduction of the $i$-signature. The removable $i$-nodes corresponding to the $-$ signs in the reduced $i$-signature are the \emph{normal} $i$-nodes of $\la$. Let $\nor i\la$ denote the number of normal $i$-nodes of $\la$, and $\renor i\la$ the partition obtained by removing all the normal $i$-nodes. Then $\renor i\la$ is a $p$-restricted partition, and we have the following.

\begin{lemma}\label{simprest}
Suppose $\la$ is a $p$-restricted partition and $i\in\zpz$. Then $\eps i{\sid\la}=\nor i\la$, and $\emax i\sid\la\cong\sid{\renor i\la}$.
\end{lemma}

\begin{eg}
Suppose $p=3$ and $\la=(6,5,3,2,1^3)$. The Young diagram of $\la$, with the residues of addable and removable nodes marked, is as follows.
\Yaddables1\Yinternals0
\[
\yngres(3,6,5,3,2,1^3)
\]
We see that the $0$-signature of $\la$ is $+--+-$. So the reduced $0$-signature is $+--$, and the normal $0$-nodes are $(2,5)$ and $(7,1)$. Hence $\eps0{\sid\la}=2$, and $\e_0^{(2)}\sid\la\cong\sid{(6,4,3,2,1^2)}$.
\end{eg}

These results will be very helpful in finding lower bounds for the number of composition factors of a Specht module $\spe\la$: since $\e_i^{(r)}$ is an exact functor, we have $\eps iT\ls\eps i{\spe\la}$ for any composition factor $T$ of $\spe\la$; furthermore, the number of composition factors $T$ (with multiplicity) for which equality holds must equal the composition length of $\spe{\reall i\la}$. Hence we have the following result.

\begin{lemma}\label{213}
Suppose $\la$ is a partition and $i\in\zpz$. Then the composition length of $\spe\la$ is at least the composition length of $\spe{\reall i\la}$, with equality if and only if $\e_i^{\remo i\la}(S)\neq0$ for every composition factor $S$ of $\spe\la$.
\end{lemma}

Analogously to the restriction functors $\e_i$, one can define \emph{induction functors} $\f_i:\bbf\sss n\mood\to\bbf\sss{n+1}\mood$ and obtain similar combinatorial results. In particular, we have the following analogue of \cref{213}, which will occasionally be helpful.

\begin{lemma}\label{213f}
Suppose $\la$ is a partition and $i\in\zpz$, and let $\addall i\la$ denote the partition obtained from $\la$ by adding all the addable $i$-nodes. Then the composition length of $\spe\la$ is at least the composition length of $\spe{\addall i\la}$.
\end{lemma}

We now generalise our notation slightly: given $i_1,\dots,i_r\in\zpz$, we write
$\reall{i_1\dots i_r}\la$ to mean $\reall{i_r}{(\dots\reall{i_2}{(\reall{i_1}\la)}\dots)}$. Clearly, \cref{213} holds with $\reall{i_1\dots i_r}\la$ in place of $\reall i\la$. If $\la$ is $p$-restricted, we write $\renor{i_1\dots i_r}\la$ to mean $\renor{i_r}{(\dots\renor{i_2}{(\renor{i_1}\la)}\dots)}$.

\subsection{Restriction and the Mullineux map}

It follows fairly easily from the definition of the $\e_i$ that they behave well with respect to the functor $-\otimes\sgn$. In fact, the following \lcnamecref{mullrest} is the basis for Kleshchev's combinatorial algorithm for computing the Mullineux map.

\begin{lemma}\label{mullrest}
Suppose $i\in\zpz$ and $M\in\bbf\sss n\mood$. Then
\[
\e_i(M\otimes\sgn)\cong(\e_{-i}M)\otimes\sgn.
\]
\end{lemma}
As a consequence of this \lcnamecref{mullrest} and \cref{simprest}, we have $\nor i\la=\nor{-i}{\mull p\la}$ for any $p$-restricted $\la$.

\subsection{Restriction and restrictisation}\label{rrsubsec}

In this section we show how to extract more information from the preceding discussion on restriction functors using \cref{jrest}. Given $i\in\zpz$, we write $\noregdown i\la$ if $\e_i^{\remo i\la}\sid{\la\rest}=0$; that is, if there is a power of $\e_i$ which kills $\sid{\la\rest}$ but not $\spe\la$. By \cref{simprest}, this is just the combinatorial condition $\nor i{\la\rest}<\remo i\la$. More generally, given residues $i_1,\dots,i_r$, we write $\noregdown{i_1\dots i_r}\la$ if 
\[
\e_{i_r}^{\remo{i_r}{\reall{i_1\dots i_{r-1}}\la}}\dots\e_{i_2}^{\remo{i_2}{\reall{i_1}\la}}\e_{i_1}^{\remo{i_1}\la}\sid{\la\rest}=0,
\]
which is the same as saying that for some $1\ls l\ls r$
\[
\nor{i_l}{\renor{i_1\dots i_{l-1}}{(\la\rest)}}<\remo{i_l}{\reall{i_1\dots i_{l-1}}{\la}}.
\]

The following is an immediate consequence of \cref{jrest,213}.

\begin{lemma}\label{ligh0}
Suppose $\la$ is a partition and that $\nrd i$ for some $i\in\zpz$. Then $\spe\la$ is reducible.
\end{lemma}

We now prove a more complicated version of this \lcnamecref{ligh0} for self-conjugate partitions. We assume for the moment that $p$ is odd, and write $p=2\pt+1$. Given $0\ls i\ls\pt$ and a \sc{} partition $\la$, we define $\al i$ to be the partition obtained by repeatedly removing all removable nodes of residue $i$ and $-i$; then $\al i$ is self-conjugate, and
\[
\al i=\begin{cases}
\reall i\la&(i=0)\\
\reall{i(-i)}\la=\reall{(-i)i}\la&(0<i<\pt)\\
\reall{i(-i)i}\la=\reall{(-i)i(-i)}\la&(i=\pt).
\end{cases}
\]

Now we have the following, which will be our main inductive tool for proving \cref{main}.

\begin{propn}\label{ligh1}
Suppose $p=2\pt+1$ is an odd prime and $\la$ is a self-conjugate partition, and that one of the following occurs.
\begin{enumerate}
\item\label{liga}
There is some $0\ls i\ls\pt$ such that $\spe{\al i}$ has at least three \cf s.
\item\label{ligb}
There is some $0\ls i<\pt$ such that $\spe{\al i}$ is irreducible.
\item\label{ligc}
There is some $0\ls i\ls \pt$ such that $\spe{\al i}$ is reducible and either $\nrd i$ or $\nrd{-i}$.
\item\label{lige}
$\nrd{\pt(-\pt)\pt}$ and $\nrd{(-\pt)\pt(-\pt)}$.
\end{enumerate}
Then the composition length of $\spe\la$ is not $2$.
\end{propn}

\begin{pfenum}
\item
This follows by applying \cref{213} one, two or three times.
\item
Assume $i>0$; a similar but simpler argument applies in the case $i=0$. Since $i<\pt$, there is no $i$-node adjacent to a $(-i)$-node, and so $\remo i{\reall{-i}\la}=\remo i\la$. So
\[
\spe{\al i}=\e_i^{(\remo i\la)}\e_{-i}^{(\remo{-i}\la)}\spe\la.
\]
If $\spe\la$ has exactly two composition factors, then by \cref{scmull} these are of the form $S$ and $S\otimes\sgn$. Since $\spe{\al i}$ is irreducible and the restriction functors are exact, exactly one of $\e_i^{\remo i\la}\e_{-i}^{\remo{-i}\la}S$ and $\e_i^{\remo i\la}\e_{-i}^{\remo{-i}\la}(S\otimes\sgn)$ must be non-zero. But by \cref{mullrest} and the fact that (since $i<h$) $\e_i$ and $\e_{-i}$ commute, we have
\[
\e_i^{\remo i\la}\e_{-i}^{\remo{-i}\la}(S\otimes\sgn)\cong
(\e_{-i}^{\remo{-i}\la}\e_i^{\remo i\la}S)\otimes\sgn=
(\e_i^{\remo i\la}\e_{-i}^{\remo{-i}\la}S)\otimes\sgn,
\]
a contradiction.
\item
Suppose $\nrd i$. Then $\e_i^{\remo i\la}\sid{\la\rest}=0$, so that by \cref{213} $\spe{\reall i\la}$ has strictly fewer composition factors than $\spe\la$. Hence $\spe{\al i}$ has strictly fewer composition factors than $\spe\la$. Since $\spe{\al i}$ is reducible, this means that $\spe\la$ has at least three \cf s.

The case where $\nrd{-i}$ is proved in the same way.
\item
If $\spe\la$ has exactly two composition factors, then by \cref{scmull,jrest} these are $\sid{\la\rest}$ and $\sid{\la\rest}\otimes\sgn$. If $\nrd{\pt(-\pt)\pt}$ and $\nrd{(-\pt)\pt(-\pt)}$, then
\[
\e_\pt^{\remo\pt{\reall{\pt(-\pt)}\la}}\e_{-\pt}^{\remo{-\pt}{\reall\pt\la}}\e_\pt^{\remo\pt\la}\sid{\la\rest}=
0=\e_{-\pt}^{\remo{-\pt}{\reall{(-\pt)\pt}\la}}\e_\pt^{\remo\pt{\reall{-\pt}\la}}\e_{-\pt}^{\remo{-\pt}\la}\sid{\la\rest},
\]
and the second equality together with \cref{mullrest} and the fact that $\la$ is \sc{} gives
\[
\e_\pt^{\remo\pt{\reall{\pt(-\pt)}\la}}\e_{-\pt}^{\remo{-\pt}{\reall\pt\la}}\e_\pt^{\remo\pt\la}(\sid{\la\rest}\otimes\sgn)=
0,
\]
so that by exactness
\[
\spe{\al\pt}=\e_\pt^{(\remo\pt{\reall{\pt(-\pt)}\la})}\e_{-\pt}^{(\remo{-\pt}{\reall\pt\la})}\e_\pt^{(\remo\pt\la)}\spe\la=
0,
\]
a contradiction.
\end{pfenum}

In order to use \cref{ligh1}, it will be very helpful to be able to test the condition $\nor i{\la\rest}<\remo i\la$ for a given partition $\la$ without having to construct $\la\rest$. We now prove a new result which will enable us to do this.

Recall the definition of the ramps used in the definition of $\la\rest$, and write
\begin{align*}
\lad l\la&\text{ for the number of nodes of $\la$ in ramp $l$,}\\
\alad l\la&\text{ for the number of addable nodes of $\la$ in ramp $l$, and}\\
\rlad l\la&\text{ for the number of removable nodes of $\la$ in ramp $l$,}
\end{align*}
setting all of these numbers to be zero when $l<0$.

\begin{lemma}\label{ladar}
For any $\la$ and any $l$, we have
\[
\alad l\la-\rlad{l-p}\la=\delta_{l0}-\lad l\la+\lad{l-1}\la+\lad{l-p+1}\la-\lad{l-p}\la.
\]
Hence if $\la$ and $\mu$ are partitions with $\la\rest=\mu\rest$, then $\alad l\la-\rlad{l-p}\la=\alad l\mu-\rlad{l-p}\mu$ for every $l$.
\end{lemma}

\begin{pf}
We assume $p>2$; a modification to the argument is required when $p=2$, and since this case is not relevant to the main results in this paper, we feel content to leave this case to the reader. We also assume $l>0$, with the case $l=0$ being trivial.

Suppose $(r,c)$ is a node in ramp $l$. Assuming first that $r,c>1$, we have nodes $(r,c-1)$, $(r-1,c)$ and $(r-1,c-1)$ in ramps $l-1$, $l-p+1$ and $l-p$ respectively. By checking the possible cases for which of these four nodes are nodes of $\la$, we easily find that the formula holds when restricted just to these four nodes. If $r=1$, then the same argument applies looking at just the two nodes $(1,c)$ and $(1,c-1)$, and a similar statement applies when $c=1$. Summing over all $(r,c)$ in ramp $l$ gives the result.

The second sentence of the lemma follows immediately, since if $\la\rest=\mu\rest$, then $\lad m\la=\lad m\mu$ for all $m$.
\end{pf}

\begin{propn}\label{renorl}
Suppose $\la$ is a partition, and $i\in\{0,\dots,p-1\}$. Then $\nor i{\la\rest}\ls\remo i\la$, with equality if and only if $\la$ does not have a removable $i$-node and an addable $i$-node in a later ramp.

Furthermore, if equality occurs then
\[
\renor i{(\la\rest)}=(\reall i\la)\rest.
\]
\end{propn}

\begin{rmksenum}
\item
Of course, the inequality $\nor i{\la\rest}\ls\remo i\la$ follows from \cref{jrest,spechtrest,simprest}. But it naturally comes out of the argument below, and it is interesting to have a purely combinatorial proof.
\item
The final statement of the \lcnamecref{renorl} (together with \cref{simprest}) shows that given a partition $\la$ and residues $i_1,\dots,i_r$, we have $\noregdown{i_1\dots i_r}\la$ if and only if for some $1\ls l\ls r$
\[
\nor{i_l}{(\reall{i_1\dots i_{l-1}}{\la})\rest}<\remo{i_l}{\reall{i_1\dots i_{l-1}}{\la}}.
\]
This will be very useful in the proof of our main theorem.
\end{rmksenum}

\begin{pf}[Proof of \cref{renorl}]
For this proof, we introduce some notation. Given a non-negative integer $m$, we let ${+}[m]$ denote a string of $m$ $+$ signs, and $-[m]$ a string of $m$ $-$ signs. Given any integer $m$, we set
\[
{\pm}[m]=\begin{cases}
{+}[m]&(m\gs0)\\
-[-m]&(m\ls0).
\end{cases}
\]
Note that when constructing the reduction of any sign sequence, we may at any point replace an interval ${-}[a]\,{+}[b]$ with ${\pm}[b-a]$.

Now let $\mu=\la\rest$. Since $\mu$ is $p$-restricted, addable and removable $i$-nodes of $\mu$ lie above addable and removable $i$-nodes in later ramps, and within a given ramp the addable nodes lie above the removable nodes. Hence the $i$-signature of $\mu$ has the form
\[
{+}[\alad i\mu]\,{-}[\rlad i\mu]\,{+}[\alad{i+p}\mu]\,{-}[\rlad{i+p}\mu]\,{+}[\alad{i+2p}\mu]\,{-}[\rlad{i+2p}\mu]\,\cdots.
\]
So the reduced $i$-signature of $\mu$ is the reduction of the sequence
\[
\scrs={+}[\alad i\mu]\,{\pm}[\alad{i+p}\mu-\rlad i\mu]\,{\pm}[\alad{i+2p}\mu-\rlad{i+p}\mu]\cdots.
\]
Now consider the sign sequence
\[
\scrl={+}[\alad i\la]\,{-}[\rlad i\la]\,{+}[\alad{i+p}\la]\,{-}[\rlad{i+p}\la]\,{+}[\alad{i+2p}\la]\,{-}[\rlad{i+2p}\la]\,\cdots.
\]
The reduction of $\scrl$ is the reduction of the sequence
\[
{+}[\alad i\la]\,{\pm}[\alad{i+p}\la-\rlad i\la]\,{\pm}[\alad{i+2p}\la-\rlad{i+p}\la]\cdots,
\]
but by \cref{ladar} the latter sequence coincides with $\scrs$. So the reduced $i$-signature of $\mu$ is just the reduction of $\scrl$, and hence the number of $-$~signs in the reduced $i$-signature is at most the number of $-$ signs in $\scrl$, which is $\remo i\la$. So $\nor i\mu\ls\remo i\la$. Equality occurs if and only if $\scrl$ is already reduced, i.e.\ every $+$ in $\scrl$ occurs before every $-$; this is the same as saying that $\la$ does not have a removable $i$-node and an addable $i$-node in a later ramp.

For the case where equality occurs, the condition on $\la$ means that we can find $L\equiv i\ppmod p$ such that $\rlad{L-ap}\la=\alad{L+ap}\la=0$ for all $a>0$. This implies in particular that
\begin{align*}
\alad{L+(a+1)p}\la-\rlad{L+ap}\la&\gs0\qquad\text{ for all $a<0$,}\\
\alad{L+(a+1)p}\la-\rlad{L+ap}\la&\ls0\qquad\text{ for all $a\gs0$.}
\end{align*}
By \cref{ladar} these inequalities also hold with $\mu$ in place of $\la$. So the reduced $i$-signature of $\mu$ has the form
\[
\begin{array}c
\hbox to \textwidth{${+}[\alad i\mu]\,{+}[\alad{i+p}\mu-\rlad i\mu]\ \cdots\ {+}[\alad L\mu-\rlad{L-p}\mu]$\hfil}\\\hfill {-}[\rlad L\mu-\alad{L+p}\mu]\ {-}[\rlad{L+p}\mu-\alad{L+2p}\mu]\ \cdots.
\end{array}
\]
We see that the $-$ signs in the reduced $i$-signature correspond to removable nodes in ladders $L,L+p,L+2p,\dots$; more precisely, $\mu$ has $\rlad{L+ap}\mu-\alad{L+(a+1)p}\mu=\rlad{L+ap}\la$ normal nodes in ramp $L+ap$ for each $a\gs0$. So $\renor i\mu$ is obtained by removing $\rlad{L+ap}\la$ nodes from ramp $L+ap$ for each $a$; hence $\renor i\mu=(\reall i\la)\rest$.
\end{pf}

\begin{eg}
Take $p=3$ and $\la=(14,5,2^3,1^5)$, so that $\mu=\la\rest=(6,5,4^2,3^2,2,1^3)$. The Young diagrams of these partitions, with the residues of nodes and addable nodes marked, are as follows.
\Yaddables1
\[
\yngres(3,14,5,2^3,1^5)\qquad\yngres(3,6,5,4^2,3^2,2,1^3)
\]
Taking $i=0$, we get the following values.
\[
\begin{array}{cc@{\ }c@{\ }c@{\ }c@{\ }c@{\ }c}\hline
l&\lad l\la&\alad l\la&\rlad l\la&\lad l\mu&\alad l\mu&\rlad l\mu\\\hline
0&1&0&0&1&0&0\\
3&2&0&0&2&0&0\\
6&3&1&1&3&1&1\\
9&2&0&1&2&0&1\\
12&2&0&0&2&0&1\\
15&0&0&0&0&1&0\\
18&1&0&1&1&0&1\\\hline
\end{array}
\]
So $\la$ satisfies the hypothesis in \cref{renorl}. We have
\[\alad{3a+3}\la-\rlad{3a}\la=\alad{3a+3}\mu-\rlad{3a}\mu\]
for all $a$, and this is non-negative for $a\ls1$ and non-positive for $a\gs2$. The $0$-signature of $\mu$, with the signs labelled according to the ramps containing the corresponding nodes, is
\[
+_6-_6-_9-_{12}+_{15}-_{18}.
\]
So $\mu$ has three normal nodes, in ramps $6$, $9$ and $18$. Since the removable $0$-nodes of $\la$ lie in ramps $6$, $9$ and $18$, we have $\renor0{(\la\rest)}=(\reall 0\la)\rest$.

\end{eg}

\section{The abacus}\label{absec}

In this section we describe the abacus notation for partitions, which we shall use in the proof of our main theorem. We also discuss \emph{\qs{} partitions}, which label Specht modules whose composition factors are well understood.

\subsection{The abacus display for a partition}\label{abdispsec}

We assume throughout this section that $p<\infty$, and fix an abacus with $p$ infinite vertical runners, which we number $0,\dots,p-1$ from left to right. We mark positions $\dots,-2,-1,0,1,2,\dots$ on the runners, reading from left to right along successive rows, with runner $i$ containing the positions congruent to $i$ modulo $p$. We say that position $m$ is \emph{later} than position $l$ (or position $l$ is \emph{earlier} than position $m$) if $m>l$. For example, if $p=5$ then the positions are marked as follows.
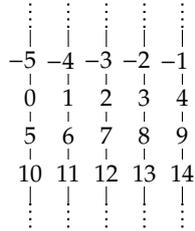
\[
\begin{tikzpicture}[xscale=.5,yscale=.5,every node/.style={fill=white,inner sep=1.8pt}]
\foreach\x in{0,...,4}{\draw[thick,dotted](\x,.5)--(\x,-.2);\draw(\x,-.2)--(\x,-4.8);\draw[thick,dotted](\x,-4.8)--(\x,-5.5);}
\draw(0,-1)node{\footnotesize$\mathllap-5$}++(1,0)node{\footnotesize$\mathllap-4$}++(1,0)node{\footnotesize$\mathllap-3$}++(1,0)node{\footnotesize$\mathllap-2$}++(1,0)node{\footnotesize$\mathllap-1$}++(-4,-1)node{\footnotesize$0$}++(1,0)node{\footnotesize$1$}++(1,0)node{\footnotesize$2$}++(1,0)node{\footnotesize$3$}++(1,0)node{\footnotesize$4$}++(-4,-1)node{\footnotesize$5$}++(1,0)node{\footnotesize$6$}++(1,0)node{\footnotesize$7$}++(1,0)node{\footnotesize$8$}++(1,0)node{\footnotesize$9$}++(-4,-1)node{\footnotesize$10$}++(1,0)node{\footnotesize$11$}++(1,0)node{\footnotesize$12$}++(1,0)node{\footnotesize$13$}++(1,0)node{\footnotesize$14$};
\end{tikzpicture}
\]
Now given a partition $\la$, we place a bead on the abacus in position $\la_i-i$ for each $i\in\bbn$. The resulting configuration is called the \emph{abacus display} for $\la$. We call a position \emph{occupied} if there is a bead at that position, and \emph{vacant} otherwise; we may also say that there is a \emph{space} in position $l$ if that position is vacant.

\begin{eg}
Take $\la=(12,10,9,7,5,4,3^2,2,1^7)$ and $p=5$. Then the abacus display for $\la$ is as follows (when drawing abacus displays, we will always suppress the numbering of positions).
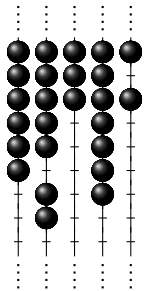
\[
\abacus(vvvvv,bbbbb,bbbbn,bbbbb,bbnbn,bbnbn,bnnbn,nbnbn,nbnnn,nnnnn;vvvvv)
\]
\end{eg}

Abacus combinatorics are well-established, so we quote some facts without further explanation. Taking an abacus display for $\la$ and sliding all the beads up their runners as far as possible, we obtain an abacus display for the $p$-core of $\la$. The $p$-weight of $\la$ is the number of pairs $l<m$ such that $l\equiv m\ppmod p$, position $l$ is vacant and position $m$ is occupied.

Now we consider $p$-quotients. Given an abacus display for a partition $\la$, we define a partition $\quo i\la$ by examining runner $i$ in isolation as a $1$-runner abacus display and reading off the corresponding partition. In other words, $\quo i\la_j$ equals the number of vacant positions above the $j$th lowest bead on runner $i$. The $p$-tuple $\pqd\la$ is called the \emph{$p$-quotient} of $\la$, and from the last paragraph the sum $|\quo0\la|+\dots+|\quo{p-1}\la|$ equals the $p$-weight of $\la$.

\begin{eg}
Taking the abacus display from the last example and sliding beads up their runners, we obtain an abacus display for the $5$-core of $\la$, namely $(9,8,6,5^2,4,3^2,2,1^2)$.
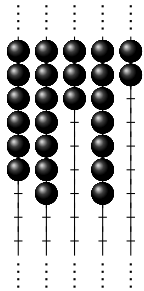
\[
\abacus(vvvvv,bbbbb,bbbbb,bbbbn,bbnbn,bbnbn,bbnbn,nbnbn,nnnnn,nnnnn;vvvvv)
\]
From the abacus display for $\la$, we see that the $5$-quotient of $\la$ is $\big(\varnothing,(1^2),\varnothing,\varnothing,(1)\big)$, so the $5$-weight of $\la$ is $3$.
\end{eg}

Since we shall be dealing with self-conjugate partitions, we record the following fact: if $\la$ is a \sc{} partition, then for any $l$ there is a bead in position $l$ in the abacus display for $\la$ if and only if there is a space in position $-l-1$. Moreover, the $p$-quotient $\pqd\la$ satisfies $\quo i\la={\quo{p-1-i}\la}'$ for each $i$.

\subsection{Addable and removable nodes and ramps}

Suppose $\la$ is a partition and $i\in\zpz$. The removable $i$-nodes of $\la$ correspond to the occupied positions $l$ on runner $i$ in the abacus display for $\la$ for which position $l-1$ is vacant. Removing such a node corresponds to moving the bead from position $l$ to position $l-1$. Similarly, addable $i$-nodes correspond to occupied positions $l-1$ on the abacus display with position $l$ vacant, and adding such a node corresponds to moving the bead from position $l-1$ to position $l$. The order of the addable and removable $i$-nodes from top to bottom of the Young diagram for $\la$ is the same as the order of the corresponding positions on runner $i$ \emph{from bottom to top}.

We record a minor lemma that we shall use later.

\begin{lemma}\label{addrem}
Suppose $\la$ is a partition, and that for some $i\in\zpz$ $\la$ has at least one addable $i$-node and at least one removable $i$-node. Then $\quo{i-1}\la\neq\varnothing$ or $\quo i\la\neq\varnothing$.
\end{lemma}

\begin{pf}
If $\quo{i-1}\la=\quo i\la=\varnothing$, then there are integers $a,b$ such that in the abacus display for $\la$:
\begin{itemize}
\item
position $kp+i-1$ is occupied if and only if $k\ls a$;
\item
position $lp+i$ is occupied if and only if $l\ls b$.
\end{itemize}
If $a\gs b$, then $\la$ has no removable $i$-nodes, while if $a\ls b$, then $\la$ has no addable $i$-nodes; contradiction.
\end{pf}

Now we use the description of addable and removable nodes on the abacus to enable us to apply \cref{renorl} using the abacus.

\begin{lemma}\label{arab}
Suppose $\la$ is a partition, $i\in\zpz$ and $k,l$ are integers with $k<l$. Suppose that positions $kp+i-1$ and $lp+i$ are occupied in the abacus display for $\la$, while positions $kp+i$ and $lp+i-1$ are vacant. Suppose furthermore that the number of beads in positions $kp+i+1,\dots,lp+i-2$ is at least $l-k$. Then \nrd i.
\end{lemma}

\begin{pf}
From the above discussion, the occupied position $lp+i$ with a space in position $lp+i-1$ corresponds to a removable $i$-node $(c,\la_c)$, while the vacant position $kp+i$ with a bead in position $kp+i-1$ corresponds to an addable $i$-node $(d,\la_d+1)$, with $c<d$. In fact, the definition of the abacus display means that $d=c+b+1$, where $b$ is the number of beads in positions $kp+i+1,\dots,lp+i-2$. The definition of the abacus display gives $\la_c-c=lp+i$ and $\la_d-d=kp+i-1$, and hence the removable node $(c,\la_c)$ lies in ramp $(c+l-1)p+i$, and the addable node $(d,\la_d+1)$ lies in ramp $(d+k-1)p+i$. The difference between these ramp numbers is
\[
(d+k)p-(c+l)p=(b+1-l+k)p
\]
which is strictly positive by assumption. Hence $\la$ has a removable $i$-node and an addable $i$-node in a later ladder, and so \nrd i by \cref{renorl}.
\end{pf}

\subsection{Rouquier partitions}\label{rouqsec}

We now describe a class of partitions for which the corresponding decomposition numbers are very well understood. In order to do this, it will be helpful to impose a new ordering on the runners of the abacus; this approach was first taken by Richards \cite{rich}.

Given a partition $\la$, construct the abacus display for the $p$-core of $\la$. Let $q_i(\la)$ be the first vacant position on runner $i$, for each $i$. We may write $q_i(\la)$ just as $q_i$ if $\la$ is understood, and we make the observation that $q_0+\dots+q_{p-1}=\binom p2$.

Let $\pi=\pi_\la$ be the unique permutation of $\{0,\dots,p-1\}$ such that $q_{\pi(0)}<\dots<q_{\pi(p-1)}$. We say that runner $\pi(0)$ is the \emph{smallest} runner in the abacus display, and runner $\pi(p-1)$ the \emph{largest}. We define $\oquo i\la=\quo{\pi(i)}\la$ for each $i$, and we define the \emph{\opq} of $\la$ to be $\opqd\la$.

\begin{eg}
Continuing from the last example with $\la=(12,10,9,7,5,4,3^2,2,1^7)$ and $p=5$, we find that $(q_0,q_1,q_2,q_3,q_4)=(5,11,-8,13,-11)$, so that $\pi$ is the $5$-cycle $(0,4,3,1,2)$, and the \opq{} of $\la$ is $\big[(1),\varnothing,\varnothing,(1^2),\varnothing\big]$.
\end{eg}

Now say that $\la$ is a \emph{Rouquier partition} if $q_{\pi(i)}-q_{\pi(i-1)}>(w-1)p$ for all $1\ls i<p$, where $w$ is the $p$-weight of $\la$. The composition factors of Specht modules labelled by Rouquier partitions are relatively well understood, thanks to the work of Chuang and Tan \cite{ct} and Turner \cite{turner}. In order to state this result, we introduce some notation. Suppose $\la$, $\mu$ and $\nu$ are Rouquier partitions with the same $p$-core and $p$-weight, and with \opq s $\opqd\la$, $\opqd\mu$ and $\opqd\nu$ respectively, and that $\oquo{p-1}\mu=\oquo{p-1}\nu=\varnothing$. Define $d_{\la\mu}$ to be the sum, over all choices of partitions $\sigma^{(0)},\dots,\sigma^{(p-2)}$ and $\tau^{(1)},\dots,\tau^{(p-1)}$, of
\[
\prod_{i=0}^{p-1}\lir{\oquo i\la}{\tau(i)'}{\sigma(i)}\prod_{i=0}^{p-2}\lir{\oquo i\mu}{\sigma(i)}{\tau(i+1)}.
\]
Here, $\lir\alpha\beta\gamma$ denotes the Littlewood--Richardson coefficient (which is to be regarded as zero if $|\alpha|\neq|\beta|+|\gamma|$) and the partitions $\tau(0)$ and $\sigma(p-1)$ should be read as $\varnothing$. Also define
\[
a_{\mu\nu}=\begin{cases}
\prod_{i=0}^{p-2}[\weyl{\oquo i\mu}\,:\,\schursimp{\oquo i\nu}]&(\text{if $|\oquo i\mu|=|\oquo i\nu|$ for all $i$})\\
0&(\text{otherwise}).
\end{cases}
\]

Now we have the following statement.

\begin{thmciting}{\textup{\textbf{\cite[Theorem 1.1]{ct}, \cite[Theorem 29]{turner}.}}\ }\label{rouqdecomp}
Suppose $\la$ and $\nu$ are Rouquier partitions with the same $p$-core and $p$-weight, and with \opq s $\opqd\la$ and $\opqd\nu$ respectively. Then:
\begin{enumerate}
\item\label{rd1}
$\nu$ is $p$-restricted if and only if $\oquo{p-1}\nu=\varnothing$;
\item
if $\nu$ is $p$-restricted, then
\[
\dn\la\nu=\sum_\mu d_{\la\mu}a_{\mu\nu},
\]
summing over all $p$-restricted partitions $\mu$ with the same $p$-core and $p$-weight as $\la$.
\end{enumerate}
\end{thmciting}

\subsection{\qs{} partitions}

In this section we generalise the notion of a Rouquier partition, to define a class of Specht modules whose composition length can be deduced from \cref{213,213f,rouqdecomp}. This material will be familiar to many experts, although it does not seem to have appeared in this form before.

Suppose $\la$ is a partition and construct the abacus display for $\la$. Say that $\la$ is \emph{\qs} if the following property holds: there do not exist runners $i\neq j$ such that the first space on runner $i$ is earlier than the last bead on runner $j$ and the first space on runner $j$ is earlier than the last bead on runner $i$. This terminology was first used by James and Mathas \cite{jmq-1} in the case $p=2$.

We assemble some facts about \qs{} partitions. Given a partition $\la$, let $q_i=q_i(\la)$ for $i\in\zpz$, and $\pi=\pi_\la$.

\begin{propn}\label{qsstuff}
Suppose $\la$ is a partition, and $i,j\in\zpz$.
\begin{enumerate}
\item\label{pt1}
If $\la$ is \qs{} and $q_i<q_j$, then the first space on runner $j$ is later than the last bead on runner $i$.
\item\label{pt2}
If $\la$ is \qs{} and has at least two addable $i$-nodes, then $q_{i-1}>q_i$.
\item\label{soisaddall}
If $\la$ is \qs{}, then so is $\addall i\la$.
\item
$\la$ is Rouquier if and only if every partition with the same $p$-core and $p$-weight as $\la$ is \qs.
\end{enumerate}
\end{propn}

\begin{pfenum}
\item
By construction, the first space on runner $i$ is in position $q_i-p\len{{\quo i\la}}$, and the last bead on runner $j$ is in position $q_j-p+p\quo j\la_1$. So either the first space on runner $i$ is earlier than the last bead on runner $j$ (in which case the first space on runner $j$ is later than the last bead on runner $i$ by the \qs{} property), or $q_j-q_i<p$ and $\quo i\la=\quo j\la=\varnothing$. But in this case the first space on runner $j$ is in position $q_j$ and the last bead on runner $i$ is in position $q_i-p$, so the required result still holds.
\item
Since $\la$ has at least two addable nodes, there are integers $k<l$ such that positions $kp+i-1$ and $lp+i-1$ are occupied while positions $kp+i$ and $lp+i$ are vacant. This implies that the first space on runner $i$ is earlier than the last bead on runner $i-1$, so $q_{i-1}>q_i$ by part (\ref{pt1}).
\item
We assume $\addall i\la$ is not \qs, and show that $\la$ is not either. Since $\addall i\la$ is not \qs, there are $k\neq l$ such that in the abacus display for $\addall i\la$ the first space on runner $k$ occurs before the last bead on runner $l$, and vice versa. Clearly if $\{k,l\}\cap\{i-1,i\}=\emptyset$, then $\la$ is not \qs, so assume otherwise.

Suppose $k=i-1$ and $l=i$, and that in the abacus display for $\addall i\la$ the first space on runner $i$ is in position $ap+i$, and the last bead on runner $i-1$ in position $bp+i-1$; then by assumption $a<b$. But now the definition of $\addall i\la$ means that in the abacus display for $\la$ positions $ap+i-1$ and $ap+i$ are both vacant, and positions $bp+i-1$ and $bp+i$ are both occupied. So the first space on runner $i-1$ is earlier than the last bead on runner $i$ and vice versa, so $\la$ is not \qs.

Next consider the case where $k=i-1$ and $l\neq i$. In the abacus display for $\addall i\la$ the first space on runner $l$ is earlier than the last bead on runner $i-1$, which is in position $bp+i-1$, say. In the abacus display for $\la$, positions $bp+i-1$ and $bp+i$ are both occupied, so the first space on runner $l$ is earlier than the last bead on runner $i-1$ and the last bead on runner $i$. Now suppose that in the abacus display for $\addall i\la$ the first space on runner $i-1$ is in position $ap+i-1$; by assumption this is earlier than the last bead on runner $l$. In the abacus display for $\la$ at least one of the positions $ap+i-1$ and $ap+i$ is vacant, so the last bead on runner $l$ is later than either the first space on runner $i-1$ or the first space on runner $i$. Either way, $\la$ fails to be \qs.

The case where $k=i$ and $l\neq i-1$ works in a very similar way.
\item
First we show that if $\la$ is Rouquier, then $\la$ is \qs. Take $i,j\in\zpz$, and suppose without loss that $q_i<q_j$. As in the proof of (\ref{pt1}), the first space on runner $j$ is in position $q_j-p\len{{\quo j\la}}$, and the last bead on runner $i$ is in position $q_i-p+p\quo i\la_1$. Now
\[
\len{{\quo j\la}}+\quo i\la_1\ls|\quo i\la|+|\quo j\la|\ls w,
\]
where $w$ is the $p$-weight of $\la$. Since $q_j-q_i>(w-1)p$, the last bead on runner $i$ is earlier than the first space on runner $j$. So $\la$ is \qs. Since every partition with the same $p$-core and $p$-weight as $\la$ is also Rouquier, it follows that every such partition is also \qs.

If $\la$ is not Rouquier, take $i,j$ such that $q_i<q_j$ and $q_j-q_i<(w-1)p$. Then from the last paragraph is it clear how to construct a partition with the same $p$-core and $p$-weight as $\la$ which is not \qs.
\end{pfenum}

Our aim is to show that the composition length of a Specht module labelled by a \qs{} partition is the same as that of a Specht module labelled by a Rouquier partition with the same \opq. The following \lcnamecref{qsindstep} gives us the inductive step.

\begin{propn}\label{qsindstep}
Suppose $\la$ is a \qs{} partition and $i\in\zpz$, with $q_{i-1}>q_i$. Then:
\begin{enumerate}
\item
$\la$ has at least one addable $i$-node and has no removable $i$-nodes;
\item
$\addall i\la$ is \qs{} and has the same \opq{} as $\la$;
\item
$\spe\la$ and $\spe{\addall i\la}$ have the same composition length.
\end{enumerate}
\end{propn}

\begin{pfenum}
\item
$\la$ cannot have a removable $i$-node, since then in the abacus display for $\la$ there would be a space on runner $i-1$ with a bead in a later position on runner $i$, contradicting \cref{qsstuff}(\ref{pt1}). If $\la$ does not have an addable $i$-node either, then for every integer $k$, position $kp+i-1$ is occupied if and only if position $kp+i$ is occupied. Clearly then the same is true in the abacus display for the $p$-core of $\la$, which gives $q_{i-1}=q_i-1$, a contradiction.
\item
The fact that $\addall i\la$ is \qs{} has already been proved in \cref{qsstuff}(\ref{soisaddall}), so we just check that $\la$ and $\addall i\la$ have the same \opq. Since $\la$ has no removable $i$-nodes, the abacus display for $\addall i\la$ is obtained from that for $\la$ by `swapping runners $i-1$ and $i$'; that is, position $kp+i$ is occupied in the abacus display for $\addall i\la$ if and only if position $kp+i-1$ is occupied in the abacus display for $\la$, and vice versa. Hence we have
\[
\quo j{(\addall i\la)}=\begin{cases}
\quo{i-1}\la&(j=i)\\
\quo i\la&(j=i-1)\\
\quo j\la&(\text{otherwise}),
\end{cases}
\qquad
q_j(\addall i\la)=\begin{cases}
q_{i-1}(\la)+1&(j=i)\\
q_i(\la)-1&(j=i-1)\\
q_j(\la)&(\text{otherwise}),
\end{cases}
\]
and the latter statement gives $\pi_{\addall i\la}=\perc(\imo i)\circ\pi_\la$. Hence $\addall i\la$ has the same \opq{} as $\la$.
\item
Since $\la$ has no removable $i$-nodes, this follows from \cref{213,213f}.
\end{pfenum}

Now we can prove the main result of this section.

\begin{propn}\label{qssamelength}
Suppose $\la$ is a \qs{} partition. Then there is a Rouquier partition $\mu$ with the same \opq{} as $\la$, and $\spe\la$ and $\spe\mu$ have the same composition length.
\end{propn}

\begin{pf}
Define
\[
d_j=\left\lfloor\frac{q_{\pi(j)}-q_{\pi(j-1)}}p\right\rfloor
\]
for $1\ls j<p$, so that $\la$ is Rouquier if $d_i\gs w-1$ for all $i$. We shall proceed by induction on
\[
M:=\sum_{j=1}^{p-1}\max\{0,w-1-d_j\},
\]
and for fixed $d_1,\dots,d_{p-1}$ we use downwards induction on $N:=\sum_iq_i^2$. Since $\sum_iq_i$ is constant, $N$ is bounded for fixed $(d_1,\dots,d_{p-1})$; so it suffices to show that we can replace $\la$ with a partition which yields a smaller value of $M$, or the same values of $d_1,\dots,d_{p-1}$ and a larger value of $N$.

Suppose that for some $i,k\in\zpz$ we have $q_{i-1}>q_k>q_i$, and consider the partition $\addall i\la$. By \cref{qsindstep}, $\addall i\la$ has the same \opq{} as $\la$, and the Specht modules $\spe\la$ and $\spe{\addall i\la}$ have the same composition length. So we can replace $\la$ with $\addall i\la$. From the proof of \cref{qsindstep} we see that making this replacement does not change any of the integers $d_j$ and strictly increases $N$.

Alternatively, assume there are no such $i,k$. The only way this can happen is if there is some $l\in\zpz$ such that $q_l>q_{l+1}>\dots>q_{l-1}$. But now if $\la$ is not Rouquier, take a $j\gs1$ such that $d_j<w-1$, and let $i=\pi(j-1)$. Then $\pi(j)=i-1$, and again we consider the partition $\addall i\la$. Again, we can replace $\la$ with $\addall i\la$ by \cref{qsindstep}, and now the replacement increases $d_j$ by $1$ while fixing all the other $d_k$, and in particular decreases $M$.
\end{pf}

\begin{eg}
Continuing from the last example, we see that $\la=(12,10,9,7,5,4,3^2,2,1^7)$ is $5$-quotient-separated, but not Rouquier, since $\la$ has $5$-weight $3$, and $(d_1,d_2,d_3,d_4)=(0,2,1,0)$. Following the proof of \cref{qssamelength}, we have $q_1>q_0>q_2$, so we can replace $\la$ with the partition $\addall2\la=(13,10^2,7,5,4^2,3,2^2,1^6)$, which has the following abacus display.
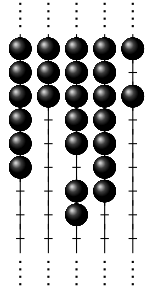
\[
\abacus(vvvvv,bbbbb,bbbbn,bbbbb,bnbbn,bnbbn,bnnbn,nnbbn,nnbnn,nnnnn;vvvvv)
\]
We now have $q_3>q_0>q_4$, so we can replace $\la$ with $\addall4\la$. We can continue in this way until we reach the partition $\la=(14^2,11^2,7,5^3,3^3,1^8)$, which has abacus display
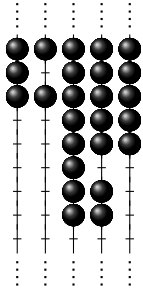
\[
\qquad\abacus(vvvvv,bbbbb,bnbbb,bbbbb,nnbbb,nnbbb,nnbnn,nnbbn,nnbbn,nnnnn;vvvvv)
\]
and $q_2>q_3>q_4>q_0>q_1$ (and still $(d_1,d_2,d_3,d_4)=(0,2,1,0)$). Now we replace $\la$ with $\addall i\la$ for $i=1,3$ or $4$, and continue.
\end{eg}

\subsection{\jm s and R-partitions on the abacus}\label{rpartabsec}

Since we shall be using the abacus extensively in the proof of our main theorem, it will be useful to have characterisations of \jm s and R-partitions in terms of their abacus displays. We take these from \cite{mfirred,altred}.

\begin{propnc}{mfirred}{Proposition 2.1}\label{jmabacus}
Suppose $\la$ is a partition with \opq{} $\opqd\la$. Then $\la$ is a \jm{} if and only if all the following hold.
\begin{itemize}
\item
$\la$ is \qs.
\item
$\oquo0\la$ is a $p$-restricted \jm.
\item
$\oquo{p-1}\la$ is a $p$-regular \jm.
\item
$\oquo i\la=\varnothing$ for $i\neq0,p-1$.
\end{itemize}
\end{propnc}

\begin{propnc}{altred}{Proposition 4.11}\label{roneabacus}
Suppose $p=2\pt+1$ is an odd prime and $\la$ is a \sc{} partition. Then $\la$ is an \rone{} if and only if there is some $j\gs0$ such that the abacus display for $\la$ satisfies the following (equivalent) conditions.
\begin{itemize}
\item
Position $jp+\pt$ is occupied, while all other positions later than position $\pt-1$ are vacant.
\item
Position $-1-jp-\pt$ is vacant, while all other positions earlier than position $-\pt$ are occupied.
\end{itemize}
\end{propnc}

\begin{propnc}{altred}{Proposition 4.12}\label{rtwoabacus}
Suppose $p=2\pt+1$ is an odd prime and $\la$ is a \sc{} partition with \opq{} $\opqd\la$. Then $\la$ is an \rtwo{} if and only if all of the following hold.
\begin{itemize}
\item
$\la$ is \qs.
\item
$\oquo0\la$ is a $p$-restricted \jm.
\item
$\oquo{p-1}\la$ is a $p$-regular \jm.
\item
$\oquo\pt\la=(1)$.
\item
$\oquo i\la=\varnothing$ for $i\neq0,\pt,p-1$.
\end{itemize}
\end{propnc}

We observe an important consequence of \cref{roneabacus} for the proof of our main theorem.

\begin{propn}\label{roneremall}
Suppose $p=2\pt+1$ is an odd prime and $0\ls i<h$, and that $\la$ is a \sc{} partition such that $\al i$ is an \rone. Then $\la$ is an \rone.
\end{propn}

\begin{pf}
In general, the abacus display for $\la$ is obtained from the abacus display for $\al i$ by moving beads from runner $i-1$ to the adjacent positions on runner $i$, and moving beads from runner $-i-1$ to the adjacent positions on runner $-i$. Since $\al i$ is an \rone, \cref{roneabacus} shows that the only possible beads that can be moved are in positions $i-1$ and $-i-1$. But even if beads are moved from these positions to positions $i$ and $-i$, then the resulting abacus display still satisfies the conditions in \cref{roneabacus}, so $\la$ is an \rone.
\end{pf}

\subsection{Proof of \cref{main} for \qs{} partitions}

We are now ready to prove \cref{main} in the case where $\la$ is \qs. Throughout this section we assume that $p$ is an odd prime, and write $p=2\pt+1$. We begin with some very simple facts about Littlewood--Richardson coefficients.

\begin{lemma}\label{lrlem}\indent
\begin{enumerate}
\vspace{-\topsep}
\item\label{lrlem1}
Suppose $\alpha$ and $\beta$ are non-empty partitions. Then there are at least two partitions $\gamma$ for which $\lir\gamma\alpha\beta>0$.
\item\label{lrlem2}
Suppose $\gamma$ is a partition with $|\gamma|>1$. Then there are at least three ordered pairs $(\alpha,\beta)$ of partitions such that $\lir\gamma\alpha\beta>0$.
\end{enumerate}
\end{lemma}

\begin{pf}
There are various `Littlewood--Richardson rules' for computing the coefficients $\lir\gamma\alpha\beta$, for example, the version given (in \cite[Theorem A1.3.3]{ec2} and elsewhere) in terms of \emph{Littlewood--Richardson tableaux} of shape $\gamma\setminus\alpha$. The results in this proof are easy to see using this rule.
\begin{enumerate}
\item
Set $\gamma=(\alpha_1+\beta_1,\alpha_2+\beta_2,\dots)$ and let $\delta$ be the partition obtained by putting the sequence $(\alpha_1,\beta_1,\alpha_2,\beta_2,\dots)$ in decreasing order. Then $\lir\gamma\alpha\beta$ and $\lir\delta\alpha\beta$ are both non-zero; the fact that $\alpha$ and $\beta$ are both non-empty ensures that $\gamma\neq\delta$; indeed, $\gamma_1=\alpha_1+\beta_1>\max\{\alpha_1,\beta_1\}=\delta_1$.
\item
$(\gamma,\varnothing)$ and $(\varnothing,\gamma)$ are obviously two such pairs. For a third, choose a removable node of $\gamma$, and let $\alpha$ be the partition obtained by removing this node. Then $\alpha\neq\varnothing$ (since $|\la|>1$) and we may take $(\alpha,(1))$ as our third pair.\qedhere
\end{enumerate}
\end{pf}

First we make an observation about the \opq{} of a \sc{} partition. Suppose $\la$ is \sc; then the $p$-core of $\la$ is also \sc, and together with the last paragraph of \cref{abdispsec} this implies that $q_i+q_{p-1-i}=p-1$ for every $i$. Hence the permutation $\pi=\pi_\la$ satisfies $\pi(p-1-i)=p-1-\pi(i)$ for every $i$. As a consequence, the \opq{} exhibits the same symmetry as the $p$-quotient, namely that $\oquo{p-1-i}\la={\oquo i\la}'$ for every~$i$.

\begin{propn}\label{mainqs}
Suppose $\la$ is a \sc{} \qs{} partition and that $\spe\la$ has exactly two composition factors. Then $\la$ is an \rtwo, and if $p=3$ then $\la$ has $3$-weight $1$.
\end{propn}

\begin{pf}
To begin with, we assume $\la$ is a Rouquier partition, and use \cref{rouqdecomp}. Let $\opqd\la$ be the \opq{} of $\la$.

For each pair of partitions $\alpha,\beta$, we fix an arbitrary choice of a partition $\pon\alpha\beta$ such that $\lir{\pon\alpha\beta}\alpha\beta>0$, and if $\alpha\neq\varnothing\neq\beta$, then (appealing to \cref{lrlem}(\ref{lrlem1})) we fix an arbitrary choice of a partition $\ptw\alpha\beta\neq\pon\alpha\beta$ such that $\lir{\ptw\alpha\beta}\alpha\beta>0$.

Suppose first that there is some $1\ls i\ls\pt-1$ such that $\oquo i\la\neq\varnothing$. Consider the three partitions $\mu^1,\mu^2,\mu^3$ with the same $p$-core as $\la$, and with \opq s
\begin{align*}
&\big[\oquo0\la,\dots,\oquo{i-2}\la,\pon{\oquo{i-1}\la}{{\oquo i\la}'},{\oquo{i+1}\la}',\dots,{\oquo{p-1}\la}',\varnothing\big],\\
&\big[\oquo0\la,\dots,\oquo{i-1}\la,\pon{\oquo i\la}{{\oquo{i+1}\la}'},{\oquo{i+1}\la}',\dots,{\oquo{p-1}\la}',\varnothing\big],\\
&\big[\oquo0\la,\dots,\oquo{p-2-i}\la,\pon{\oquo{p-1-i}\la}{{\oquo{p-i}\la}'},{\oquo{p+1-i}\la}',\dots,{\oquo{p-1}\la}',\varnothing\big].
\end{align*}
Since these \opq s are distinct, the partitions $\mu^1,\mu^2,\mu^3$ are distinct. By \cref{rouqdecomp}(\ref{rd1}), $\mu^1,\mu^2,\mu^3$ are $p$-restricted, with $d_{\la\mu^k}>0$, and hence (since $a_{\mu\mu}=1$ for any $\mu$) $\dn\la{\mu^k}>0$ for $k=1,2,3$. So $\spe\la$ has at least three composition factors; contradiction.

So we have $\oquo i\la=\varnothing$ for $1\ls i\ls\pt-1$, and hence also for $\pt+1\ls i\ls p-2$. Now suppose that $\oquo\pt\la=\varnothing$. Then the only $p$-restricted partition $\mu$ for which $d_{\la\mu}>0$ is the partition with the same $p$-core as $\la$ and with \opq{}
\[
\big[\oquo0\la,\varnothing,\dots,\varnothing,{\oquo{p-1}\la}',\varnothing\big],
\]
for which we have $d_{\la\mu}=1$. Since $\oquo0\la={\oquo{p-1}\la}'$, the sum $\sum_\nu a_{\mu\nu}$ is a perfect square (namely, the square of the composition length of the Weyl module $\weyl{\oquo0\la}$). So the composition length of $\spe\la$ is a perfect square; contradiction.

So we have $\oquo\pt\la\neq\varnothing$. Now we look at the case $p=3$, and observe that in this case $\oquo0\la=\oquo2\la=\varnothing$; for if not, then the partitions $\mu$ with the same $p$-core as $\la$ and with \opq s
\begin{align*}
&\big[\pon{\oquo0\la}{{\oquo1\la}'},{\oquo2\la}',\varnothing\big],\\
&\big[\ptw{\oquo0\la}{{\oquo1\la}'},{\oquo2\la}',\varnothing\big],\\
&\big[\oquo0\la,\pon{\oquo1\la}{{\oquo2\la}'},\varnothing\big]
\end{align*}
all give $d_{\la\mu}>0$ and hence $\dn\la\mu>0$.

Next we suppose (with $p$ arbitrary again) that $|\oquo\pt\la|>1$. Then by \cref{lrlem}(\ref{lrlem2}) we can find three pairs of partitions $\alpha,\beta$ such that $\lir{\oquo\pt\la}\alpha\beta>0$. For each such pair, we construct the partition $\mu$ with the same $p$-core as $\la$ and with \opq
\[
\big[\oquo0\la,\varnothing,\dots,\varnothing,\alpha',\beta,\varnothing,\dots,\varnothing,{\oquo{p-1}\la}',\varnothing\big]
\]
(where the $\alpha',\beta$ occur in positions $\pt-1,\pt$), and we have $d_{\la\mu}>0$ and hence $\dn\la\mu>0$. So again $\spe\la$ has at least three composition factors, a contradiction.

So we have $\oquo\pt\la=(1)$. This completes the analysis in the case $p=3$. In the case $p\gs5$, the only $p$-restricted $\mu$ for which $d_{\la\mu}>0$ are those with \opq
\[
\big[\oquo0\la,\varnothing,\dots,\varnothing,(1),\varnothing,\dots,\varnothing,{\oquo{p-1}\la}',\varnothing\big],
\]
where the $(1)$ occurs either in position $\pt-1$ or in position $\pt$. For each of these two partitions we have $d_{\la\mu}=1$, while $\sum_\nu a_{\mu\nu}$ is again equal to the square of the composition length of $\weyl{\oquo0\la}$. So $\weyl{\oquo0\la}$ must be simple, i.e.{} $\oquo0\la$ is \ $p$-restricted \jm. So by \cref{rtwoabacus} $\la$ is an \rtwo.
\end{pf}

\section{Homomorphisms between Specht modules}\label{homsec}

In this section, we review some results on homomorphisms between Specht modules and prove a result that we shall need later. This material is discussed at length elsewhere, so in the interests of brevity we specialise as much as possible.

\subsection{Tableau homomorphisms}\label{tabhomsubsec}

We begin with some combinatorics. Throughout this section let $\la$ and $\mu$ be fixed compositions of $n$. A \emph{$\la$-tableau of type $\mu$} is a function $T$ from the Young diagram of $\la$ to $\bbn$ with the property that exactly $\mu_i$ nodes are mapped to $i$, for each $i$. We write $T_{r,c}$ for the image of the node $(r,c)$ under $T$, and we illustrate $T$ by drawing the Young diagram and filling the $(r,c)$-box with $T_{r,c}$, for each $(r,c)$. A tableau is \emph{row-standard} if its entries weakly increase from left to right along the rows, and \emph{semistandard} if it is row-standard and its entries strictly increase down the columns.

Recall that $\yper\la$ denotes the Young permutation module associated with $\la$. For each row-standard $\la$-tableau $T$ of type $\mu$, there is an $\bbf\sss n$-homomorphism $\tht T:\yper\la\to\yper\mu$ defined in \cite[\S13]{jbook}. If $\la$ is a partition, then the Specht module $\spe\la$ is a submodule of $\yper\la$, and the restriction of $\tht T$ to $\spe\la$ is denoted $\hth T$. If $T$ is semistandard, we refer to $\hth T$ as a \emph{semistandard homomorphism}. Now we have the following result.

\begin{thmc}{jbook}{Lemma 13.11 \& Theorem 13.13}\label{homback}
Suppose $\la$ and $\mu$ are compositions of $n$.
\begin{enumerate}
\item
The set
\[
\lset{\tht T}{T\text{ a row-standard $\la$-tableau of type $\mu$}}
\]
is an $\bbf$-basis for $\ho{\yper\la}{\yper\mu}$.
\item
If $\la$ is a partition, the set
\[
\lset{\hth T}{T\text{ a semistandard $\la$-tableau of type $\mu$}}
\]
is an $\bbf$-basis for the space of all $\bbf\sss n$-homomorphisms $\spe\la\to\yper\mu$ which can be extended to $\yper\la$.
\item
If $\la$ is a partition and $p\gs3$, then every $\bbf\sss n$-homomorphism $\spe\la\to\yper\mu$ can be extended to $\yper\la$.
\end{enumerate}
\end{thmc}
In view of this \lcnamecref{homback}, a natural way to construct a homomorphism $\spe\la\to\spe\mu$ when $p\gs3$ is to find a linear combination of semistandard homomorphisms $\spe\la\to\yper\mu$ whose image lies in $\spe\mu$; the latter condition can be checked using James's Kernel Intersection Theorem \cite[Corollary 17.18]{jbook} and the author's results from \cite{mfgarnir}, so that we now have a reasonably fast algorithm \cite{mfgarnir} for computing $\ho{\spe\la}{\spe\mu}$. Even when $p=2$, this method can often be used to construct homomorphisms between Specht modules (including the homomorphisms described in this paper), though it will not in general find all homomorphisms.

In this section we want to show the existence of a non-zero homomorphism $\spe\la\to\spe\mu$ in a certain case, and we construct this as the composition of two known homomorphisms between Specht modules. But we have some work to do in showing that this composition is non-zero. In order to do this, we need to discuss \emph{dominance}. If $T$ is a row-standard $\la$-tableau of type $\mu$, let $\rec Tlr$ denote the total number of entries less than or equal to $l$ in rows $1,\dots,r$ of $T$. If $U$ is another row-standard $\la$-tableau of type $\mu$, we say that $T$ \emph{dominates} $U$ (and write $T\dom U$) if $\rec Tlr\gs \rec Ulr$ for all $l,r$.

\begin{eg}
The dominance order on the set of row-standard $(3,2)$-tableaux of type $(2^2,1)$ may be represented by the following Hasse diagram.
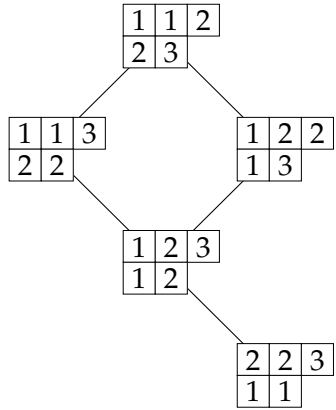
\newlength\xpos\newlength\ypos
\[
\begin{tikzpicture}
\foreach\x in {1.5}{
\draw(0.5,0.5)--++(-2*\x,2*\x)--++(\x,\x)--++(\x,-\x)--++(-\x,-\x);
\setlength\xpos{0cm}\setlength\ypos{0.5cm}
\tyoung(\xpos,\ypos,223,11)
\addtolength\xpos{-\x cm}\addtolength\ypos{\x cm}
\tyoung(\xpos,\ypos,123,12)
\addtolength\xpos{\x cm}\addtolength\ypos{\x cm}
\tyoung(\xpos,\ypos,122,13)
\addtolength\xpos{-\x cm}\addtolength\ypos{\x cm}
\tyoung(\xpos,\ypos,112,23)
\addtolength\xpos{-\x cm}\addtolength\ypos{-\x cm}
\tyoung(\xpos,\ypos,113,22)}
\end{tikzpicture}
\]
\end{eg}

The following \lcnamecref{domhom} will be very useful.

\begin{lemma}\label{domhom}
Suppose $T$ is a row-standard $\la$-tableau, and write $\hth T$ as a linear combination $\sum_St_S\hth S$ of semistandard homomorphisms. Then $S\dom T$ for each $S$ with $t_S\neq0$.
\end{lemma}

\begin{pf}
This follows directly from the algorithm given in \cite[\S5.2]{mfgarnir} for semistandardising $\hth T$ (and in fact is quite easy to see from \cite[\S13]{jbook} where the homomorphisms $\hth T$ are introduced).
\end{pf}

\begin{eg}
Continuing from the last example, if we let $T=\young(122,13)$, then by \cref{domhom} $\hth T$ should be a linear combination of homomorphisms $\hth S$ with $S\dom T$. From the diagram above we see that the only such $S$ is $S=\young(112,23)$. And indeed (as can easily be shown using the results below) $\hth T=-\hth S$.
\end{eg}

Now we describe two particular constructions of homomorphisms that we shall use. Throughout this section we assume that $p$ is finite.

\subsection{One-node Carter-Payne homomorphisms}\label{cp}

Suppose $\la$ is a partition of $n$ with a removable node $(a,b)$ and an addable node $(c,d)$ of the same residue, with $c>a$. Let $\mu$ be the partition obtained by removing $(a,b)$ and adding $(c,d)$. Then there is a non-zero $\bbf\sss n$-homomorphism $\spe\la\to\spe\mu$; this is a special case of the Carter--Payne Theorem \cite{cp}, and an explicit formula for this homomorphism may be found in the paper of Lyle \cite{slcp}. For simplicity, we concentrate on a special case.

Suppose $\la$ and $\mu$ are as above, and suppose additionally that $\la$ has no removable nodes in rows $a+1,\dots,c-1$; that is, $\la_{a+1}=\dots=\la_{c-1}=d$. For each $a<r\ls c$ define a $\la$-tableau $\cp\la\mu r$ of type $\mu$ by
\[
\cp\la\mu r_{x,y}=\begin{cases}
r&(\text{if $(x,y)=(a,b)$})\\
x+1&(\text{if $r\ls x<c$ and $y=d$})\\
x&(\text{otherwise}).
\end{cases}
\]
\begin{eg}
Taking $p=3$, $\la=(4^2,2^3,1^2)$, $\mu=(4,3,2^4,1)$, $(a,b)=(2,4)$ and $(c,d)=(6,2)$, we have
\[
\cp\la\mu3=\young(1111,2223,34,45,56,6,7),\qquad
\cp\la\mu4=\young(1111,2224,33,45,56,6,7),\qquad
\cp\la\mu5=\young(1111,2225,33,44,56,6,7),\qquad
\cp\la\mu6=\young(1111,2226,33,44,55,6,7).
\]
\end{eg}

Now we have the following result, which is a very special case of \cite[Theorem 4.5.4]{slcp}.

\begin{propn}\label{cpprop}
Suppose $\la$, $\mu$ are as above. Then $\ho{\spe\la}{\spe\mu}$ is one-dimensional, and a non-zero homomorphism $\spe\la\to\spe\mu$ is given by $\sum_{r=a+1}^c(-1)^r\hth{\cp\la\mu r}$.
\end{propn}

\subsection{Restrictisation homomorphisms}\label{magicsubsec}

In this section we consider the homomorphisms arising in the following \lcnamecref{reghom}, which may be regarded as a homomorphism-space analogue of \cref{jrest}.

\begin{thmc}{flm}{Theorem 1.5}\label{reghom}
Suppose $\la$ is a partition of $n$. Then
\[
\dim_\bbf\left(\ho{\spe\la}{\spe{\la\rest}}\right)=1.
\]
\end{thmc}

A non-zero homomorphism $\spe\la\to\spe{\la\rest}$ is constructed in \cite{flm} as a tableau homomorphism $\hth T$, but for a tableau $T$ which is not necessarily semistandard. This leads to an additional problem (which creates a large part of the work in \cite{flm}) of showing that $\hth T$ is non-zero. In order to prove our result on composition of homomorphisms, we shall find the least dominant tableau occurring when $\hth T$ is expressed as a linear combination of semistandard homomorphisms; we note that this gives a new proof that $\hth T$ is non-zero.

First we must describe the tableau $T$. In fact, there is a range of possibilities for $T$, yielding homomorphisms $\hth T$ which agree up to sign. These tableaux are called \emph{magic tableaux} in \cite{flm}, and we define them here using one of the recursive characterisations given in \cite{flm}, which we re-phrase for our own purposes.

Let $\la$ be a partition of $n$. Write $\full\la$ for the number of `full' ramps in $\la$ (i.e.\ ramps in which every node is a node of $\la$); then it is easy to see that $\full\la=\la\rest_1$. If we look at the first ramp which is not full, we can find a node of this ramp, in row $m$ say, which is not a node of $\la$. Then we have $\la_m+(m-1)(p-1)=\full\la$, while for any $1\ls l\ls m$ we have $\la_l+(l-1)(p-1)\gs\full\la$. Call such a value of $m$ a \emph{nice} value for $\la$. Having chosen a nice value $m$, we define a partition $\la\stre$ by
\[
\la\stre_i=\begin{cases}
\la_i-p+1&(i<m)\\
\la_{i+1}&(i\gs m).
\end{cases}
\]
It is easy to see that $\left(\la\stre\right)\rest=(\la\rest_2,\la\rest_3,\dots)$; in particular, $\full{\la\stre}=\la\rest_2$. Now given any $\la\stre$-tableau $U$, define a $\la$-tableau $\insm U$ by
\[
\insm U_{x,y}=\begin{cases}
1&(\text{if $x<m$ and $y<p$})\\
U_{x,y-p+1}+1&(\text{if $x<m$ and $y\gs p$})\\
1&(\text{if $x=m$})\\
U_{x-1,y}+1&(\text{if $x>m$}).
\end{cases}
\]
Now we can define magic tableaux recursively: the unique tableau for the empty partition is magic, and if $\la\neq\varnothing$, a magic tableau for $\la$ is any tableau of the form $\insm U$, where $U$ is a magic $\la\stre$-tableau for some nice value $m$.

\begin{eg}
Take $p=3$ and $\la=(8,6,2,1^2)$, giving $\la\rest=(6,5,4,2,1)$ and $\full\la=6$. The only nice value for $\la$ is $m=3$, giving $\la\stre=(6,4,1,1)$. One can show recursively that
\begin{align*}
&\young(112222,1133,1,4)
\\
\intertext{is a magic $\la\stre$-tableau, giving the magic $\la$-tableau}
&\young(11223333,112244,11,2,5).
\end{align*}
\end{eg}

Now we have the following \lcnamecref{magichom}.

\begin{thmc}{flm}{Theorem 2.2 \& Lemma 4.2}\label{magichom}
Suppose $\la$ is a partition, and $T$ is a magic $\la$-tableau. Then $\hth T$ defines a non-zero homomorphism from $\spe\la$ to $\spe{\la\rest}$. If $U$ is any other magic $\la$-tableau, then $\hth U=\pm\hth T$.
\end{thmc}

Now we consider expressing a `magic homomorphism' as a linear combination of semistandard homomorphisms. Given a partition $\la$, define a $\la$-tableau $\rtab\la$ as follows: let $\rtab\la_{x,y}$ equal $x$ plus the number of nodes below $(x,y)$ in the same ramp which are not nodes of $\la$. Informally, we construct $\rtab\la$ by filling each box with the number of the row that box moves to when we construct $\la\rest$ from $\la$. Hence $\rtab\la$ has type $\la\rest$.

\begin{eg}
Take $\la=(8,6,2,1^2)$ and $p=3$. Then
\[
\rtab\la=\young(11111123,222234,33,4,5).
\]
\end{eg}

Our aim in this section is to prove the following statement.

\begin{propn}\label{resttab}
Suppose $\la$ is a partition and $T$ a magic $\la$-tableau, and write $\hth T$ as a linear combination $\sum_Sa_S\hth S$ of semistandard homomorphisms. Then $a_{\rtab\la}=\pm1$, and for any $S$ with $a_S\neq0$ we have $S\dom\rtab\la$.
\end{propn}

We begin with a lemma about $\rtab\la$.

\begin{lemma}\label{rtablem}
Suppose $\la$ is a partition, and $(a,b)$ and $(c,d)$ are nodes of $\la$ with $(p-1)a+b\ls(p-1)c+d$.
\begin{enumerate}
\item
If $a=c$, then $\rtab\la_{a,b}\ls\rtab\la_{c,d}$.
\item
If $a<c$, then $\rtab\la_{a,b}<\rtab\la_{c,d}$.
\end{enumerate}
In particular, $\rtab\la$ is semistandard.
\end{lemma}

\begin{pf}
For this proof, say that a node is \emph{missing} if it not a node of $\la$.
\begin{enumerate}
\item
The conditions imply that $b\ls d$. For every missing node $(x,y)$ in the same ramp as $(a,b)$ there is a missing node $(x,y+d-b)$ in the same ramp as $(a,d)$. Hence the number of missing nodes below $(a,d)$ in the same ramp is at least the number of missing nodes below $(a,b)$ in the same ramp.
\item
Arguing as in the previous case, the number of missing nodes below row $a$ and in the same ramp as $(a,b)$ is at most the number of missing nodes below row $a$ and in the same ramp as $(c,d)$. Of the latter nodes, at most $c-a-1$ lie between rows $a$ and $c$; so the number of missing nodes below row $c$ and in the same ramp as $(c,d)$ is strictly greater than the number of missing nodes below row $a$ and in the same ramp as $(a,b)$ plus $a-c$, and this gives the result.\qedhere
\end{enumerate}
\end{pf}

Now we describe the relations in \cite{fm,mfgarnir} used to `semistandardise' homomorphisms. For these and subsequent results, we need some notation for multisets of positive integers, which we collect here.
\begin{itemize}
\item
Given a multiset $X$, let $X_i$ denote the multiplicity of $i$ as an element of $X$.
\item
Given two multisets $X,Y$, let $X\sqcup Y$ denote the multiset with $(X\sqcup Y)_i=X_i+Y_i$ for all $i$.
\item
Given a multiset $X$, let $X+1$ denote the multiset obtained from $X$ by increasing each element by $1$.
\item
Given any $l,n\in\bbn$ let $\muti ln$ denote the multiset with $n$ elements all equal to $l$.
\item
Given a tableau $T$, let $T^i$ denote the multiset of entries in the $i$th row of $T$.
\end{itemize}
With this notation in place, we can state a useful result for manipulating tableau homomorphisms. This is a combination of \cite[Theorem 3.1]{mfgarnir} and \cite[Lemma 4]{fm}.

\begin{propn}\label{maingarn}
Suppose $\la$ is a partition, $A$ is a row-standard $\la$-tableau and $1\ls h<k$. Suppose $R$, $S$, $T$ are multisets of positive integers with $R\sqcup S\sqcup T=A^h\sqcup A^k$ and $|S|>\la_h$. Let $\cals$ be the set of all pairs $(U,V)$ of multisets such that $S=U\sqcup V$ and $|R|+|U|=\la_h$. For each $(U,V)\in\cals$, let $A[U,V]$ denote the row-standard $\la$-tableau with
\[
A[U,V]^i=\begin{cases}
R\sqcup U&(\text{if }i=h)\\
T\sqcup V&(\text{if }i=k)\\
A^i&(\text{otherwise}).
\end{cases}
\]
Then
\[
\sum_{(U,V)\in\cals}\prod_{i\gs1}\mbinom{R_i+U_i}{R_i}\mbinom{T_i+V_i}{T_i}\hth{A[U,V]}=0.
\]
\end{propn}

\begin{rmk}
It is shown in \cite[\S5.2]{mfgarnir} that the relations obtained from \cref{maingarn} are sufficient to express a tableau homomorphism $\hth T$ as a linear combination of semistandard homomorphisms; in fact, only the case $k=h+1$ is required. Since the semistandard homomorphisms are linearly independent, this means that any linear relation between tableau homomorphisms is a consequence of the relations obtained from \cref{maingarn}.
\end{rmk}

Before we proceed, we note a simple corollary.

\begin{cory}\label{toomany}
Suppose $\la$ is a partition, $A$ is a row-standard $\la$-tableau and $1\ls h<k$. Suppose that for some $l\in\bbn$ we have $A^h_l+A^k_l>\la_h$. Then $\hth A=0$.
\end{cory}

\begin{pf}
Take $R=\lset{i\in A^h}{i\neq l}$, $S=\muti l{A^h_l+A^k_l}$ and $T=\lset{i\in A^k}{i\neq l}$ in \cref{maingarn}. Then $\cals=\{(U,V)\}$, where $U=\{l\}^{\la_h-|R|}$ and $V=\{l\}^{\la_k-|T|}$. So the sum in \cref{maingarn} has only one term, in which the binomial coefficients are all $1$, and we get $\hth A=0$.
\end{pf}

Now we show that the operation $U\mapsto\insm U$ preserves relations between tableau homomorphisms.

\begin{propn}\label{indstep}
Suppose $\la$ is a partition and $m$ is a nice value for $\la$, and let $\la\stre$ be as above. Suppose $\calu$ is a set of row-standard $\la\stre$-tableaux of type $\left(\la\stre\right)\rest$, and that coefficients $(c_U)_{U\in\calu}$ are chosen such that $\sum_{U\in\calu}c_U\hth U=0$. Then
\[
\sum_{U\in\calu}c_U\hth{\insm U}=0.
\]
\end{propn}

\begin{pf}
From the remark following \cref{maingarn}, it suffices to consider only relations of the form given in that \lcnamecref{maingarn}. So suppose $A,h,k,R,S,T$ are as in \cref{maingarn} (with $\la\stre$ in place of $\la$). Then we have
\begin{align*}
&\sum_{(U,V)\in\cals}\prod_{i\gs1}\mbinom{R_i+U_i}{R_i}\mbinom{T_i+V_i}{T_i}\hth{A[U,V]}=0,\\
\intertext{and we want to show}
&\sum_{(U,V)\in\cals}\prod_{i\gs1}\mbinom{R_i+U_i}{R_i}\mbinom{T_i+V_i}{T_i}\hth{\insm{A[U,V]}}=0.
\end{align*}
We consider three cases.

\begin{description}
\addtocounter{casy}1
\item[Case \arabic{casy}: $k<m$.]\indent\\\setcounter{caseone}{\value{casy}}
In this case we have $(\insm A)^h=A^h+1\sqcup\opm$ and $(\insm A)^k=A^k+1\sqcup\opm$. We define
\begin{align*}
\hat R&=R+1,\\
\hat S&=S+1\sqcup\opm,\\
\hat T&=T+1\sqcup\opm.
\end{align*}
Then $\hat R,\hat S,\hat T$ satisfy the hypotheses of \cref{maingarn} (with $\insm A$ in place of $A$), so (with the obvious definition of $\hat\cals$) we have
\begin{align*}
0&=\sum_{(\hat U,\hat V)\in\hat\cals}\prod_{i\gs1}\mbinom{\hat R_i+\hat U_i}{\hat R_i}\mbinom{\hat T_i+\hat V_i}{\hat T_i}\hth{\insm A[\hat U,\hat V]}\\
&=\sum_{(\hat U,\hat V)\in\hat\cals}\mbinom{\hat U_1}{0}\mbinom{p-1+\hat V_1}{p-1}\prod_{i\gs2}\mbinom{\hat R_i+\hat U_i}{\hat R_i}\mbinom{\hat T_i+\hat V_i}{\hat T_i}\hth{\insm A[\hat U,\hat V]}.
\end{align*}
Any term with $\hat V_1>0$ can be neglected, since $\binom{p-1+v}{p-1}\equiv0\ppmod p$ for any $0<v<p$. When $\hat V_1=0$, we have $\hat U_1=p-1$, and so
\[
(\hat U,\hat V)=(U+1\sqcup\opm,V+1)
\]
for some $(U,V)\in\cals$, and $\insm A[\hat U,\hat V]=\insm{A[U,V]}$. Conversely, if $(U,V)\in\cals$, then $(U+1\sqcup\opm,V+1)\in\hat\cals$. Hence
\begin{align*}
0&=\sum_{(U,V)\in\cals}\mbinom{p-1}{0}\mbinom{p-1}{p-1}\prod_{i\gs2}\mbinom{\hat R_i+U_{i-1}}{\hat R_i}\mbinom{\hat T_i+V_{i-1}}{\hat T_i}\hth{\insm{A[U,V]}}\\
&=\sum_{(U,V)\in\cals}\prod_{i\gs2}\mbinom{R_{i-1}+U_{i-1}}{R_{i-1}}\mbinom{T_{i-1}+V_{i-1}}{T_{i-1}}\hth{\insm{A[U,V]}}\\
&=\sum_{(U,V)\in\cals}\prod_{i\gs1}\mbinom{R_i+U_i}{R_i}\mbinom{T_i+V_i}{T_i}\hth{\insm{A[U,V]}}
\end{align*}
as required.
\addtocounter{casy}1
\item[Case \arabic{casy}: $h<m\ls k$.]\indent\\
In this case $(\insm A)^h=A^h+1\sqcup\opm$, while $(\insm A)^{k+1}=A^k+1$. We define
\begin{align*}
\hat R&=R+1,\\
\hat S&=S+1\sqcup\opm,\\
\hat T&=T+1.
\end{align*}
Now $\hat R, \hat S, \hat T$ satisfy the hypotheses of \cref{maingarn} (with $\insm A$ in place of $A$, and $k+1$ in place of $k$); proceeding as in Case \arabic{caseone}, we have
\[
0=\sum_{(\hat U,\hat V)\in\hat\cals}\mbinom{\hat U_1}{0}\mbinom{\hat V_1}{0}\prod_{i\gs2}\mbinom{\hat R_i+\hat U_i}{\hat R_i}\mbinom{\hat T_i+\hat V_i}{\hat T_i}\hth{\insm A[\hat U,\hat V]}.
\]
For any pair $(\hat U,\hat V)$ with $\hat V_1>0$ we have $\hth{\insm A[\hat U,\hat V]}=0$ by \cref{toomany} (with $l=1$ and with $m$ and $k+1$ in place of $h$ and $k$). So now we need only consider pairs $(U,V)$ with $\hat U_1=p-1$ and $\hat V_1=0$, and we can proceed as in Case \arabic{caseone}.
\addtocounter{casy}1
\item[Case \arabic{casy}: $m\ls h$.]\indent\\
In this case we have $(\insm A)^{h+1}=A^h+1$ and $(\insm A)^{k+1}=A^k+1$, and applying \cref{maingarn} with $R+1,S+1,T+1,h+1,k+1,\insm A$ in place of $R,S,T,h,k,A$ yields the result.
\qedhere
\end{description}
\end{pf}

Next we introduce a result which gives relations between tableau homomorphisms which allow us to move all the $1$s in a tableau up to the top row. We use the following variation on \cref{maingarn}; this was actually proved before \cref{maingarn}, although it follows from the latter fairly easily by induction.

\begin{propnc}{fm}{Lemma 7}\label{lemma7}
Suppose $\la$ is a partition, $B$ is a row-standard $\la$-tableau and $h,r\in\bbn$. Let $\calv$ denote the set of submultisets $V$ of $B^h$ such that $r\notin V$ and $|V|=B^{h+1}_r$. For each $V\in\calv$, let $B[V]$ denote the row-standard $\la$-tableau with
\[
B[V]^i=\begin{cases}
B^h\setminus V\sqcup\muti r{B^{h+1}_r}&(\text{if }i=h)\\
B^{h+1}\setminus\muti r{B^{h+1}_r}\sqcup V&(\text{if }i=h+1)\\
B^i&(\text{otherwise}).
\end{cases}
\]
Then
\[
\hth B=(-1)^{B^{h+1}_r}\sum_{V\in\calv}\prod_{i\gs1}\mbinom{B^{h+1}_i+V_i}{V_i}\hth{B[V]}.
\]
\end{propnc}

Using this, we can prove the following.

\begin{propn}\label{move1s}
Suppose $\la$ is a partition and $m$ is a nice value for $\la$, and define $\la\stre$ as above.
\begin{enumerate}
\item
$\hth{\insm{\rtab{\la\stre}}}$ can be expressed as a linear combination $\pm\hth{\rtab\la}+\sum_Sc_S\hth S$, where each $S$ is a row-standard $\la$-tableau with $S\doms\rtab\la$.
\item
If $T$ is a semistandard $\la\stre$-tableau with $T\doms\rtab{\la\stre}$, then $\hth{\insm T}$ can be expressed as a linear combination $\sum_Sc_S\hth S$, where each $S$ is a row-standard $\la$-tableau with $S\doms\rtab\la$.
\end{enumerate}
\end{propn}

\begin{pfenum}
\item
Let $R=\rtab{\la\stre}$. We use \cref{lemma7} to re-write $\hth{\insm R}$ by moving all the $1$s up to the top row. We apply \cref{lemma7} $m-1$ times, each time with $r=1$, taking $h=m-1,m-2,\dots,1$ in turn. At a given step, we move all the $1$s from row $h+1$ to row $h$, and move a multiset of entries greater than $1$ from row $h$ to row $h+1$. Since $R$ is semistandard, the entries in row $h$ of $R$ are all at least $h$; hence the entries not equal to $1$ in row $h$ of $\insm R$ are all at least $h+1$. Furthermore, there are at least $\full{\la\stre}-(h-1)(p-1)$ $h$s in row $h$ of $R$. We have
\begin{align*}
\full{\la\stre}
&=\la\rest_2\\
&\gs\la\rest_1-p+1\tag*{(since $\la\rest$ is $p$-restricted)}\\
&=\full\la-p+1,
\end{align*}
so we see that the first $\full\la-h(p-1)$ entries in row $h$ of $\insm R$ are equal to $h+1$. So each time we apply \cref{lemma7}, one of the terms we obtain involves moving only entries equal to $h+1$ down to row $h+1$. Taking this term at every stage we obtain the tableau $\rtab\la$, and the coefficient of $\hth{\rtab\la}$ obtained is $\pm1$; indeed, the binomial coefficients in \cref{lemma7} are always trivial, since all the entries in row $h+1$ of $\insm R$ (other than the $1$s) are strictly greater than $h+1$.

Any other term we obtain from our repeated applications of \cref{lemma7} involves moving all the $1$s up to row $1$, and moving $\full\la-h(p-1)$ entries greater than or equal to $h+1$ down from row $h$ to row $h+1$ for each $1\ls h<m$, with a strict inequality at some point. Hence the resulting tableau will strictly dominate $\rtab\la$.
\item
This case is similar to the previous one: when we apply \cref{lemma7} repeatedly, we move all the $1$s in $T$ up to row $1$, and move $\full\la-h(p-1)$ entries greater than or equal to $h+1$ down from row $h$ to row $h+1$. Since $T\doms\insm{\rtab{\la\stre}}$, any tableau resulting from this process will strictly dominate $\rtab\la$.\qedhere
\end{pfenum}

\begin{eg}
Take $p=3$, $\la=(8,6,2,1^2)$ and $m=3$. Then
\[
\rtab{\la\stre}=\young(111112,2223,3,4),\qquad\insm{\rtab{\la\stre}}=\young(11222223,113334,11,4,5).
\]
Applying \cref{lemma7} twice to move the $1$s up to the top row, we get the following possibilities.
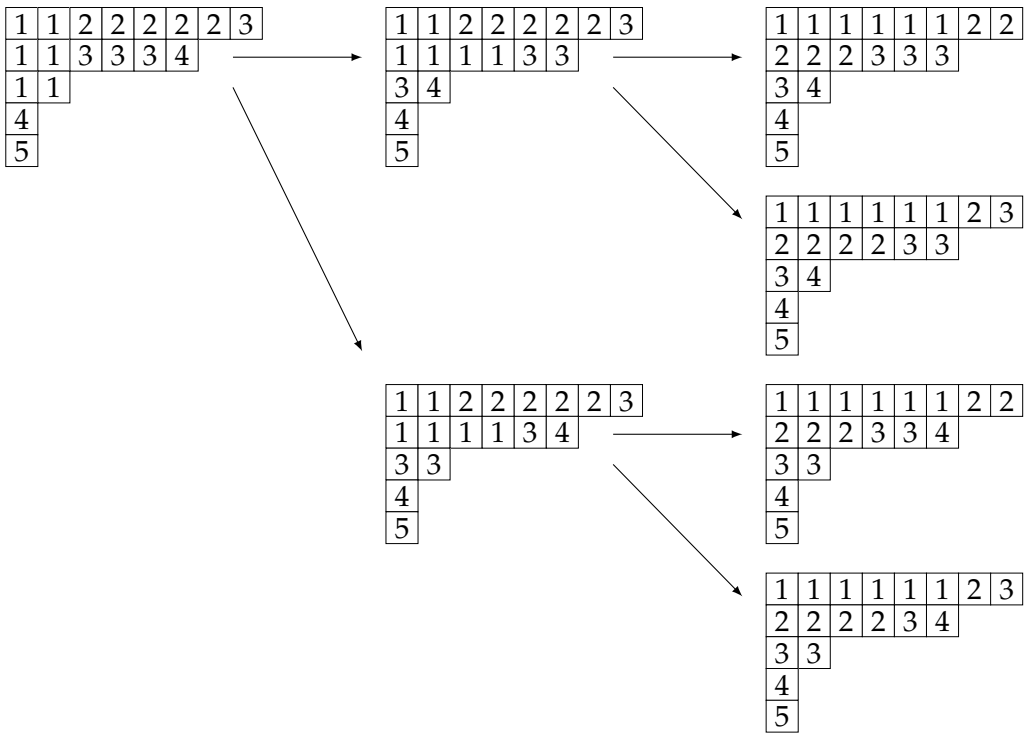
\[
\begin{tikzpicture}[scale=1,>=latex]
\foreach\x in{5}\foreach\y in{2.5}\foreach\t in {1.3}\foreach\u in{.4}\foreach\v in{1.7}{
\draw(0,0)node{\young(11222223,113334,11,4,5)};
\draw(\x,0)node{\young(11222223,111133,34,4,5)};
\draw(\x,-2*\y)node{\young(11222223,111134,33,4,5)};
\draw(2*\x,0)node{\young(11111122,222333,34,4,5)};
\draw(2*\x,-\y)node{\young(11111123,222233,34,4,5)};
\draw(2*\x,-2*\y)node{\young(11111122,222334,33,4,5)};
\draw(2*\x,-3*\y)node{\young(11111123,222234,33,4,5)};
\draw[->](\t,\u)--++(\v,0);
\draw[->](\x+\t,\u)--++(\v,0);
\draw[->](\x+\t,\u-2*\y)--++(\v,0);
\draw[->](\t,0)--++(\v,-1.4*\y);
\draw[->](\x+\t,0)--++(\v,-.7*\y);
\draw[->](\x+\t,-2*\y)--++(\v,-.7*\y);
}
\end{tikzpicture}
\]
So our homomorphism $\spe\la\to\spe{\la\rest}$ is a linear combination of the homomorphisms labelled by the semistandard tableaux at the right. The bottom tableau at the right-hand side is $\rtab\la$, occurs with coefficient $1$, and is the least dominant.
\end{eg}

\begin{pf}[Proof of \cref{resttab}]
We use induction on $|\la|$, with the case $\la=\varnothing$ being trivial. Suppose $\la\neq\varnothing$, and $m$ be the nice value chosen in the construction of $T$. Then $T=\insm U$ for a magic $\la\stre$-tableau $U$, and by induction we can assume that when we write $\hth U$ as a linear combination $\sum_Vu_V\hth V$ of semistandard homomorphisms, we have $u_{\rtab{\la\stre}}=\pm1$ while $u_V=0$ for any $V\ndom\rtab{\la\stre}$.

By \cref{indstep} we have
\[
\hth T=\hth{\insm U}=\sum_Vu_V\hth{\insm V};
\]
by \cref{move1s} this equals $\pm\hth{\rtab\la}$ plus a linear combination of homomorphisms $\hth S$ for $S\doms\rtab\la$. By \cref{domhom} each $\hth S$ can be written as a linear combination of semistandard homomorphisms $\hth R$ for $R\dom S\doms T$, and the result follows.
\end{pf}

\subsection{Composition of homomorphisms}\label{composubsec}

Our aim in this section is to show that the composition of the homomorphisms from the two previous sections is non-zero, given a certain additional condition. Specifically, we prove the following result.

\begin{propn}\label{mainhom}
Suppose $\la$ is a partition of $n$ with a removable node $(a,b)$ and an addable node $(c,d)$ of the same residue with $c>a$ and $(p-1)c+d\gs(p-1)a+b$, and that $\la$ has no removable nodes in rows $a+1,\dots,c-1$. Let $\mu$ be the partition obtained from $\la$ by removing $(a,b)$ and adding $(c,d)$. Then there is a non-zero $\bbf\sss n$-homomorphism $\spe\la\to\spe{\mu\rest}$.
\end{propn}

In order to prove \cref{mainhom}, we need to describe how to compose tableau homomorphisms. Recall that if $S$ is a tableau, then $S^j$ denotes the multiset of entries in row $j$ of $S$, and in particular $S^j_i$ denotes the number of entries equal to $i$ in row $j$ of $S$. If $x_1,x_2,\dots$ are non-negative integers with finite sum $x$, we write $(x_1,x_2,\dots)!$ for the multinomial coefficient $\mfrac{x!}{x_1!x_2!\dots}$.

\begin{propnc}{df}{Proposition 4.7}\label{tabcomp}
Suppose $\la,\mu,\nu$ are compositions of $n$, $S$ is a $\la$-tableau of type $\mu$ and $T$ is a $\mu$-tableau of type $\nu$. Let $\calx$ be the set of all collections $X=(X^{ij})_{i,j\gs1}$ of multisets such that
\[
|X^{ij}|=S^j_i\quad\text{ for each $i$, $j$,}\qquad\bigsqcup_{j\gs1}X^{ij}=T^i\quad\text{ for each $i$.}
\]
For $X\in\calx$, let $U_X$ denote the row-standard $\la$-tableau with $(U_X)^j=\bigsqcup_{i\gs1}X^{ij}$. Then
\[
\Theta_T\circ\Theta_S = \sum_{X\in\calx}\prod_{i,j\gs1}\left(X^{1j}_i,X^{2j}_i,\dots\right)!\Theta_{U_X}.
\]
\end{propnc}

Our aim is to use \cref{tabcomp} to show that the composition of the homomorphisms $\spe\la\to\spe\mu\to\spe{\mu\rest}$ from the last two sections is non-zero. Rather than attempting to give an explicit expression for this composition, we use \cref{cpprop,resttab} to find the least dominant tableau occurring when the composition is expressed in terms of semistandard homomorphisms.

Given $\la$ and $\mu$ as in \cref{mainhom}, let $V$ be the $\la$-tableau obtained from $\rtab\mu$ by moving the $(c,d)$-entry up to position $(a,b)$. Note that by \cref{rtablem} $V$ is semistandard.

\begin{lemma}\label{vtablem}
Suppose $\la$ and $\mu$ are as above, and $T$ is a semistandard $\mu$-tableau of type $\mu\rest$ with $T\doms\rtab\mu$. If $W$ is a row-standard tableau obtained from $T$ by moving an entry from row $c$ up to row $a$, then $W\ndomsby V$.
\end{lemma}

\begin{pf}
Suppose to the contrary that $V\doms W$. Since $V$ agrees with $\rtab\mu$ and $W$ agrees with $T$ in rows $1,\dots,a-1$, and we have $V\doms W$ while $\rtab\mu\domsby T$, all four tableaux must agree on rows $1,\dots,a-1$. Similarly, all four tableaux agree on rows $c+1,c+2,\dots$.

$V$ is obtained from $\rtab\mu$ by moving an entry $t=\rtab\mu_{c,d}$ from row $c$ up to row $a$. By \cref{rtablem}, $t$ is the largest entry in rows $a,\dots,c$ of $\rtab\mu$. By the previous paragraph (and since $T$ and $\rtab\mu$ have the same type) the largest entry in rows $a,\dots,c$ of $T$ is also $t$. So $W$ is obtained from $T$ by moving an entry less than or equal to $t$ from row $c$ up to row $a$. Since $T\doms\rtab\mu$, this gives $W\doms V$, a contradiction.
\end{pf}

For the next \lcnamecref{composework}, recall the tableaux $\cp\la\mu{a+1},\dots,\cp\la\mu c$ from \cref{cp}.

\begin{propn}\label{composework}
Suppose $a<r\ls c$ and $T$ is a semistandard $\mu$-tableau of type $\mu\rest$ with $T\dom\rtab\mu$.
\begin{enumerate}
\item
If $r=c$ and $T=\rtab\mu$, then $\tht T\circ\tht{\cp\la\mu r}$ equals a linear combination $\tht V+\sum_Uc_U\tht UU$, where each $U$ is a row-standard tableau with $V\ndom U$.
\item
If $r<c$ or $T\doms\rtab\mu$, then $\tht T\circ\tht{\cp\la\mu r}$ equals a linear combination $\sum_Uc_U\tht UU$, where each $U$ is a row-standard tableau with $V\ndom U$.
\end{enumerate}
\end{propn}

\begin{pf}
Taking $S=\cp\la\mu r$ and $\nu=\mu\rest$, \cref{tabcomp} simplifies considerably in our situation. It says that $\tht T\circ\tht S$ is a linear combination of row-standard tableaux $U$ obtained from $T$ by moving an entry from row $r$ to row $a$ and moving an entry from row $j+1$ to row $j$ for each $r\ls j<c$. If we let $U(T,r)$ be the particular tableau obtained by moving the largest possible entry at each stage, then we have $U\dom U(T,r)$ for any other such $U$. So it suffices to consider only the tableaux $U(T,r)$.

Since $\mu_{a+1}=\dots=\mu_{c-1}$ and $T$ is semistandard, the largest entries in rows $a+1,\dots,c-1$ are the entries in column $d$, and these are strictly increasing. Hence for $a<r<c$ we have $U(T,r)\doms U(T,r+1)$. So it suffices to consider only the tableau $U(T,c)$. But by \cref{vtablem} $V\ndom U(T,c)$ if $T\doms\rtab\mu$.

It remains to observe that the coefficient of $V$ in $\tht{\rtab\mu}\circ\tht{\cp\la\mu c}$ is $1$; but this follows from \cref{rtablem}: since $(c,d)$ lies in the same ramp as or a later ramp than $(a,b)$, we have $\rtab\mu_{c,d}>\rtab\mu_{a,b-1}$, so all the multinomial coefficients in the coefficient of $\tht V$ equal $1$.
\end{pf}

\begin{pf}[Proof of \cref{mainhom}]
By \cref{cpprop}, a non-zero homomorphism $\alpha:\spe\la\to\spe\mu$ is given by $\sum_{r=a+1}^c(-1)^r\hth{\cp\la\mu r}$. By \cref{resttab}, there is a non-zero homomorphism $\beta:\spe\mu\to\spe{\mu\rest}$ of the form $\hth{\rtab\mu}+\sum_Sb_S\hth S$, where each $S$ is a semistandard tableau strictly dominating $\rtab\mu$.

The composition $\beta\circ\alpha$ may be computed using \cref{tabcomp}. By \cref{composework}, we get
\[
\beta\circ\alpha=\pm\hth V+\sum_Uc_U\hth U,
\]
where each $U$ is a row-standard tableau with $V\ndom U$. Hence when we write such a $U$ as a linear combination of semistandard homomorphisms, the coefficient of $\hth V$ is zero, by \cref{domhom}. So when we write $\beta\circ\alpha$ as a linear combination of semistandard homomorphisms, the coefficient of $\hth V$ is $\pm1$, and in particular $\beta\circ\alpha\neq0$.
\end{pf}

\begin{rmk}
For simplicity, we have concentrated on quite a special case. It is possible to weaken the assumptions on $\la$ and $\mu$ and use the same argument: we can allow removable nodes in $\la$ between rows $a$ and $c$, as long as these none of these removable nodes has the same residue as $(a,b)$, and as long as $\mu_i-\mu_{i+1}<p$ for all $a<i<c$. We leave the reader to check the details (referring to \cite{slcp} for the formula for a homomorphism $\spe\la\to\spe\mu$).

However, it is not generally the case that the composition of a one-node Carter--Payne homomorphism and a restrictisation homomorphism is non-zero. For example, take $p=3$, $\la=(6)$ and $\mu=(5,1)$, so that $\mu\rest=(3,2,1)$. $\spe{(6)}$ is isomorphic to the simple module $\sid{(2^3)}$ (which happens to be the trivial $\bbf\sss6$-module); since $\soc{\spe{(3,2,1)}}$ is a different simple module $\sid{(3,2,1)}$, there is no non-zero homomorphism $\spe{(6)}\to\spe{(3,2,1)}$.
\end{rmk}

\section{Proof of the main theorem}\label{pfsec}

We now come to the proof of our main theorem. We proceed by induction, with our main tool being \cref{ligh1}. As we shall see, this deals with all cases except for two families of partitions which we deal with using the other techniques described above.

\subsection{Notation and assumptions}\label{assumpsec}

Throughout this section, $p=2h+1$ is an odd prime. \jmset{} denotes the set of $p$-JM-partitions, and \twofac{} denotes the set of self-conjugate partitions that (according to \cref{main}) label Specht modules with exactly two composition factors. So if $p\gs5$ then \twofac{} is the set of R-partitions, while if $p=3$ then \twofac{} comprises the \rones, the self-conjugate partitions of $3$-weight $1$ and the partition $(3^3)$.

It will be helpful to make some assumptions that will remain in force for the next few sections. Essentially, these say that $\la$ is a partition which cannot be dealt with by \cref{ligh1}(\ref{liga},\ref{ligb}) or \cref{mainqs}. Recall that if $\la$ is a partition and $0\ls i\ls\pt$, then $\al i$ denotes the partition obtained by repeatedly removing all removable nodes of residue $i$ or $-i$. Similarly, we define $\addall{\pm i}\la$ to be the partition obtained from $\la$ by repeatedly adding addable nodes of residue $\pm i$.

\smallskip
\begin{mdframed}[innerleftmargin=3pt,innerrightmargin=3pt,innertopmargin=3pt,innerbottommargin=3pt,roundcorner=5pt,innermargin=-3pt,outermargin=-3pt]
\noindent\textbf{Assumptions and notation in force for \cref{assumpsec,0casesec,icasesec,ptcase1sec,ptcase2sec}}

$\la$ is a \sc{} partition which is not \qs{} and is not in \twofac{} or \jmset. $\al i\in\twofac\cup\{\la\}$ for each $0<i<\pt$, while $\al\pt\in\twofac\cup\jmset\cup\{\la\}$.
\end{mdframed}

We observe some immediate consequences of these assumptions. Suppose $0\ls i<\pt$ and $\al i\neq\la$. Then by \cref{roneremall} \al i is not an \rone, since $\la$ is not. Also, when $p=3$ \al i cannot equal $(3^3)$, since this partition has one removable node and two addable nodes, all of residue $0$. So \al i must be an \rtwo, and if in addition $p=3$ then \al0 has $3$-weight $1$; in particular, \al i is \qs.

$\la$ is obtained from \al i by adding equal numbers of addable $i$- and $(-i)$-nodes. But we cannot have $\la=\addall{\pm i}{(\al i)}$, since then by \cref{qsstuff}(\ref{soisaddall}) $\la$ would be \qs, contradicting our assumptions. So $\la$ is obtained from $\al i$ by adding some but not all of the addable $i$-nodes, and some but not all of the addable $(-i)$-nodes. In particular, \al i has at least two addable $i$-nodes.

A similar discussion applies statement applies when $\al\pt\neq\la$. In this case $\al\pt$ can be an \rone, but only of $p$-weight $1$, since an \rone{} with $p$-weight greater than $1$ has a removable $\pt$-node. An \rone{} with $p$-weight $1$ is also an \rtwo, and so in fact we can assume \al\pt{} is an \rtwo{} or a \jm, and in particular we can assume \al\pt{} is \qs{}; as in the last paragraph, this means that $\la$ must have at least one addable $\pt$-node and at least one addable $(-\pt)$-node.

\subsection{First case}\label{0casesec}

In this section we consider the case where $\reall0\la\neq\la$.

\begin{propn}\label{0nowt}
Suppose $\al0\neq\la$ and $\quo0{(\al0)}=\varnothing$. Then \nrd0.
\end{propn}

\begin{pf}
From the discussion in \cref{assumpsec} \al0 is an \rtwo, so we have $\quo\pt{(\al0)}=(1)$, and hence $\quo\pt\la=(1)$. Since \al0 is \sc, we have $\quo{p-1}{(\al0)}=\quo0{(\al0)}=\varnothing$; so there is some $r$ such that position $kp-1$ is occupied in the abacus display for \al0 if and only if $k\ls r$, while position $kp$ is occupied if and only if $k<-r$. The fact that \al0 has at least two addable $0$-nodes means that $r\gs1$.

The abacus display for $\la$ is obtained by moving some but not all of the beads in positions $rp-1,\dots,-rp-1$ to the adjacent positions on runner $0$; since $\la$ is \sc, the bead in position $ip-1$ is moved if and only if the bead in position $-ip-1$ is moved. This means that there must be some $-r\ls i<j\ls r$ with $i\ls0$ such that $\la$ has beads in positions $jp$ and $ip-1$ and spaces in positions $jp-1$ and $ip$. If $j\gs1$, then the number of beads in positions $ip+1,\dots,jp-2$ is at least $j-i$, since for each $i<k<j$ there is a bead in position $kp$ or $kp-1$, and there is also a bead in position $\pt$; hence by \cref{arab} we have $\nrd0$. If $i\ls-2$ then again the number of intervening beads is at least $j-i$, since there is a bead in position $\pt-2p$. In the case where $j=0$, $i=-1$ and $p\gs5$, there must be a bead in position $\pt-p-1$ or $\pt-p+1$, since if both of these positions are vacant, then (by self-conjugacy) then positions $\pt-1$ and $\pt+1$ are occupied, so that \al0 is not \qs, a contradiction. So again we have at least $j-i$ beads between positions $ip$ and $jp-1$, and again $\nrd0$.

The remaining case is where $p=3$, and the only pair $i<j$ such that positions $3j$ and $3i-1$ are occupied while positions $3j-1$ and $3i$ are vacant is $(-1,0)$. In this case we must have $r=1$, and hence $\la$ has abacus display
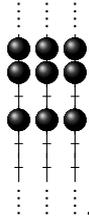
\[
\abacus(vvv,bbb,bbb,nnn,bbb,nnn,nnn;vvv).
\]
So $\la=(3^3)$, contrary to assumption.
\end{pf}

To help us deal with the case where $\quo0{(\al0)}\neq\varnothing$, we use the following lemma.

\begin{lemma}\label{0caselem}
Suppose $\al0\neq\la$. Then $\al i=\la$ for all $1<i\ls\pt$.
\end{lemma}

\begin{pf}
Since $\la$ has both addable and removable $0$-nodes, we have either $\quo0\la\neq\varnothing$ or $\quo{p-1}\la\neq\varnothing$ by \cref{addrem}. But $\la$ is self-conjugate, so in fact $\quo0\la\neq\varnothing\neq\quo{p-1}\la$.

Suppose the lemma is false, and take $1<i\ls\pt$ such that $\al i\neq\la$. Then either $\al i$ is an \rtwo, or $i=\pt$ and $\al i$ is a \jm. In either case, runner $0$ must be either the largest or smallest runner of $\al i$ by \cref{jmabacus,rtwoabacus}. But since $\la$ (and hence \al i) has an addable $0$-node, runner $0$ cannot be the smallest runner, by \cref{qsstuff}(\ref{pt1}). So runner $0$ is the largest runner of \al i, and hence (in the abacus displays for both \al i and $\la$) the first space on runner $0$ occurs after the last bead on any other runner. But \al0 is an \rtwo, so the abacus display for \al0 (and hence for $\la$) has a bead in position $\pt$; so the first space on runner $0$ in the abacus display for $\la$ is at position $p$ or later. Since $\la$ is \sc, this means that the last bead on runner $p-1$ is no later than position $-1-p$. But now $\la$ has no addable $0$-nodes, contrary to assumption.
\end{pf}

\begin{propn}\label{0somewt}
Suppose $\al0\neq\la$ and $\quo0{(\al0)}\neq\varnothing$. Then $p\gs5$, and either $\nrd0$ or $\nrd{p-1}$.
\end{propn}

\begin{pf}
The fact that $p\gs5$ follows from our standing assumptions (in particular, the definition of $\twofac$ when $p=3$). \al0 is an \rtwo{} with $\quo0{(\al0)}\neq\varnothing$, and hence $\quo\pt{(\al0)}=(1)$ and $\quo i{(\al0)}=\varnothing$ for $i\neq0,\pt,p-1$; the same applies to the $p$-quotient of $\la$. By \cref{rtwoabacus} runners $0$ and $p-1$ must be the smallest and largest runners of \al0 in some order, and in fact by \cref{qsstuff}(\ref{pt2}) runner $0$ is the smallest (since $\al0$ has at least two addable $0$-nodes). So in the abacus display for $\al0$, every position before the first space on runner $0$ is occupied, and every position after the last bead on runner $p-1$ is vacant.

$\al0$ is a self-conjugate partition, which means that position $ip-1$ in the abacus display is occupied if and only if position $-ip$ is vacant. Together with the statements in the last paragraph, this implies that there is a finite set $I$ of positive integers such that:
\begin{itemize}
\item
the occupied positions on runner $p-1$ of \al0 are the positions $ip-1$ for all $i\in I$ and all $i\ls0$; and
\item
the vacant positions on runner $0$ of \al0 are the positions $-ip$ for all $i\in I$ and all $i\ls0$.
\end{itemize}
The fact that $\al0$ has at least two addable $0$-nodes means that $I$ is non-empty.

$\la$ is obtained from \al0 by moving some but not all of the beads in positions $ip-1$ (for $i\in\pm I\cup\{0\}$) to the right; the fact that $\la$ is \sc{} means that the bead in position $ip-1$ is moved if and only if the bead in position $-ip-1$ is moved.

Suppose first that for some $i<j\ls0$ the bead in position $jp-1$ is moved but the bead in position $ip-1$ is not. Then $\la$ has beads in positions $jp$ and $ip-1$, and spaces in positions $jp-1$ and $ip$, and we claim that the number of beads in between these two positions is at least $j-i$. Indeed, there is a bead in position $kp-1$ or $kp$ (or possibly both) for each $i<k<j$, and additionally there is a bead in position $jp-\pt$ or $jp-\pt-2$ (since $\la$ is \sc{} and $\quo{\pt-1}\la=\quo{\pt+1}\la=\varnothing$). So \nrd0, by \cref{arab}.

So suppose there are no such $i$ and $j$. Then in particular the bead in position $-1$ is not moved in constructing $\la$ from \al0, but the beads in positions $\pm lp-1$ are moved, where $l=\max I$. By \cref{0caselem} $\al i=\la$ for $i\neq 0,1$, and since position $\pt-p$ is vacant, this means that positions $\pt-p+1,\pt-p+2,\dots,-2$ are vacant. On the other hand, the first paragraph of this proof shows that position $-lp-2$ must be occupied (in \al0, and hence in $\la$). So $\la$ has beads in positions $-1$ and $-lp-2$, and spaces in positions $-2$ and $-lp-1$; the number of intervening beads is at least $l$, since for each $1\ls k\ls l$ at east one of positions $-kp-1,-kp$ is occupied, and so $\nrd{p-1}$.
\end{pf}

\subsection{Second case}\label{icasesec}

Here we consider the case where there is $0<i<\pt$ such that $\al i\neq\la$ while $\al j=\la$ for all $0\ls j<i$.

\begin{propn}\label{biggestorsmallest}
Suppose $0<i<\pt$, and that $\al i\neq\la$. Then either \nrd i or \nrd{p-i}, or runner $i-1$ is the largest runner in the abacus display for \al i, or runner $i$ is the smallest runner in the abacus display for \al i.
\end{propn}

\begin{pf}
Suppose first that $\quo{i-1}{(\al i)}\neq\varnothing$. Since \al i is an \rtwo, this means that runner $i-1$ is either the largest or the smallest runner in \al i; but \al i has at least two addable $i$-nodes, so by \cref{qsstuff}(\ref{pt2}) runner $i-1$ must be the largest runner. Similarly, if $\quo i{(\al i)}\neq\varnothing$, then runner $i$ is the smallest runner in \al i.

So assume $\quo{i-1}{(\al i)}=\quo i{(\al i)}=\varnothing$. This means that there are integers $b<c$ such that in the abacus display for $\la$:
\begin{itemize}
\item
for $k<b$, both of the positions $kp+i-1$ and $kp+i$ are occupied;
\item
for $b\ls k<c$, exactly one of the positions $kp+i-1$ and $kp+i$ is occupied;
\item
for $c\ls k$, neither of the positions $kp+i-1$ and $kp+i$ is occupied.
\end{itemize}
Suppose that for some $k$, positions $kp+i$ and $(k-1)p+i-1$ are occupied while positions $kp+i-1$ and $(k-1)p+i$ are vacant. If there is at least one bead in a position between $(k-1)p+i$ and $kp+i-1$, then \nrd i, by \cref{arab}. On the other hand, if all of these positions are vacant, then in the abacus display for \al i every runner except runner $i-1$ has a space before the bead in position $kp+i-1$; so by \cref{qsstuff}(\ref{pt1}) runner $i-1$ is the largest runner of \al i.

Assume instead that there is no such $k$. Since $\la$ is obtained from \al i by adding some but not all of the addable $i$-nodes (and some but not all of the addable $(-i)$-nodes), there must be $k$ such that in the abacus display for $\la$ positions $kp+i-1$ and $(k-1)p+i$ are occupied while positions $kp+i$ and $(k-1)p+i-1$ are vacant. Since $\la$ is \sc, this means that positions $-kp+p-i+1$ and $-(k+1)p+p-i$ are occupied while positions $-kp+p-i$ and $(-k-1)p+p-i+1$ are vacant. Repeating the argument from the last paragraph (with $i$ replaced by $p-i$ and $k$ replaced by $-k$) we find that either \nrd{p-i} or runner $p-i$ is the largest runner in \al i, i.e.\ runner $i$ is the smallest.
\end{pf}

We now consider the possibilities in \cref{biggestorsmallest} in more detail. We need to treat the case $i>1$ and $i=1$ separately.

\begin{propn}\label{iwti-1}
Suppose $1<i<\pt$, and that $\al i\neq\la$ while $\al j=\la$ for all $0\ls j<i$. Suppose furthermore that runner $i-1$ is the largest runner in the abacus display for \al i. Then $\nrd i$.
\end{propn}

\begin{pf}
Since $\al i$ is an \rtwo{} and runner $i-1$ is its largest runner, we have $\quo j{(\al i)}=\varnothing$ for $j\neq i-1,\pt,p-i$. In particular, $\quo0{(\al i)}=\varnothing$. Since $\la$ (and hence $\al i$) has no removable $0$-nodes, position $0$ in the abacus is vacant (for both partitions), and hence position $lp$ is vacant for all $l\gs0$. The assumption that $\al j=\la$ for all $j<i$ then means that position $lp+i-1$ is vacant (in $\la$) for all $l\gs0$.

Let $dp+i-1$ be the first vacant position on runner $i-1$ in $\al i$; since $\al i$ is \qs{} and there is a bead in position $\pt$, we must have $d\gs1$. The fact that $dp+i-1$ is the first vacant position on runner $i-1$ means that for each $c<d$ at least one of the positions $cp+i-1$ and $cp+i$ is occupied in $\la$; in fact, since $\quo i{(\al i)}=\varnothing$, there is some $b\ls d$ such that
\begin{itemize}
\item
for $c=b,\dots,d-1$ exactly one of the positions $cp+i-1$ and $cp+i$ is occupied, and
\item
for $c<b$ both of the positions $cp+i-1$ and $cp+i$ are occupied.
\end{itemize}
By assumption $\la$ has at least one addable $i$-node, so there is some $k$ such that position $kp+i-1$ is occupied while position $kp+i$ is vacant. From the first paragraph of the proof we must have $k<0$, and in particular $k<d-1$. Taking a maximal such $k$, we therefore have positions $kp+i-1$ and $(k+1)p+i$ occupied, while positions $kp+i$ and $(k+1)p+i-1$ are vacant. There is at least one bead among positions $kp+i+1,\dots,(k+1)p+i-2$: if $k=-1$ then we can take the bead in position $-i$, while if $k<-1$ we can take the bead in position $\pt-2p$. So $\nrd i$ by \cref{arab}.
\end{pf}

Now we consider the case $i=1$.

\begin{propn}\label{1wt0}
Suppose $p\gs5$, and that $\al1\neq\la$ while $\al0=\la$. Suppose furthermore that runner $0$ is the largest runner in the abacus display for \al1. Then either $\nrd1$ or $\nrd{p-1}$.
\end{propn}

\begin{pf}
Since $\al0=\la$, position $0$ in the abacus display for $\la$ is vacant. If we let $dp$ be the first vacant position on runner $0$ in \al1, then, arguing as in the proof of \cref{iwti-1}, we must have $d\gs1$. In particular, position $0$ is occupied in the abacus display for $\al1$, and hence position $1$ is occupied in the abacus display for $\la$. Since $\la$ is \sc, position $-1$ is occupied and position $-2$ is vacant.

$\la$ has an addable $1$-node, so there is some $k$ such that position $kp$ is occupied and position $kp+1$ is vacant. If $k<0$, then taking a maximal negative such $k$ and copying the last part of the proof of \cref{iwti-1}, we get \nrd1. So suppose $k>0$. Then by self-conjugacy position $-kp-2$ is occupied and position $-kp-1$ is vacant. For every $0>l>-k$ at least one of the positions $lp-2$ and $lp-1$ is occupied (since otherwise we would have $\quo{p-2}{(\al1)}\neq\varnothing$); furthermore, at least one of positions $-p$ and $1-p$ is occupied (since $d>0$); hence the number of beads between positions $-kp-1$ and $-2$ is at least $k$, and so \nrd{p-1}.
\end{pf}

Now we consider the situation where runner $i$ is the smallest in the abacus display for $\al i$. Again, we have to consider the cases $i>1$ and $i=1$ separately.

\begin{propn}\label{iwti}
Suppose $1<i<\pt$, and that $\al i\neq\la$ while $\al j=\la$ for all $0\ls j<i$. Suppose furthermore that runner $i$ is the smallest runner in the abacus display for \al i. Then $\nrd i$.
\end{propn}

\begin{pf}
The fact that $\al i$ is an \rtwo{} with smallest runner $i$ means that $\quo j{(\al i)}=\varnothing$ for all $j\neq i,\pt,p-i-1$. Suppose that for some $l<k$ positions $kp+i$ and $lp+i-1$ are occupied in the abacus display for $\la$, while positions $kp+i-1$ and $lp+i$ are vacant; choose such a pair $(k,l)$ with $k-l$ as small as possible. Arguing as in the first paragraph of the proof of \cref{iwti-1}, we must have $l<0$. Now for every $k>m>l$, positions $mp+i$ and $mp+i-1$ must both be occupied: if both were vacant then we would have $\quo{i-1}{(\al i)}\neq\varnothing$, while if only one were vacant then the choice of $k,l$ would be contradicted. So if $k-l>1$ then the number of beads between positions $lp+i$ and $kp+i-1$ is at least $k-l$, and so $\nrd i$, as required. If $k-l=1$ and $k<0$, then there is a bead in position $lp+\pt$, so again $\nrd i$. Finally if $k=0$ and $l=-1$ then there is a bead in position $-i$, so again $\nrd i$.

So we can assume that there are no such $k$ and $l$. This means that the last bead on runner $i$ in the abacus display for $\la$ is either earlier than or adjacent to the last bead on runner $i-1$, and hence occurs no later than position $-p+i$. Furthermore, since $\la$ has at least one addable $i$-node and there are no beads on runner $i-1$ after position $-p+i-1$, there must be a space on runner $i$ no later than position $-p+i$. Now we claim that for every $i<j<\pt$, there is a space on runner $j$ no later than position $-p+j$. If this claim is not true, take the first $j$ for which it fails. Then there is a space on runner $j-1$ with an adjacent bead on runner $j$, so $\al j\neq\la$. Hence $\la$ has at least one addable $j$-node; since there is no space on runner $j$ before position $j$, there must be a bead on runner $j-1$ in position $j-1$ or later. This means that $j>i+1$ and that $\quo{j-1}\la\neq\varnothing$, and hence $\quo{j-1}{(\al i)}\neq\varnothing$, contradicting the first sentence of the proof.

So our claim holds. In particular, if $i<\pt-1$ there is a space on runner $\pt-1$ no later than position $-p+\pt-1$. So (because $\quo{\pt-1}\la=\varnothing$) there is no bead on runner $\pt-1$ after position $-2p+\pt-1$. But there is a bead in position $\pt$ and a bead in position $cp+\pt$ for all $c\ls-2$, so $\la$ has removable $\pt$-nodes but no addable $\pt$-nodes; contradiction. The same argument applies when $i=\pt-1$ and there is a space in position $-p+\pt-1$.

The only remaining case is where $i=\pt-1$ and there is a bead in position $-p+\pt-1$ and a space on runner $\pt-1$ earlier than this bead. From our assumptions so far, there is also a bead in position $-p+\pt-2$, and a space on runner $\pt-2$ earlier than this bead. Hence $\quo{\pt-2}\la\neq\varnothing$. Now we consider the partition $\al\pt$. Since $\quo{\pt-2}{(\al\pt)}=\quo{\pt-2}\la\neq\varnothing$, $\al\pt$ is either a \jm{} or an \rtwo, and runner $\pt-2$ is either the largest or the smallest runner of $\al\pt$. The first space on runner $\pt-2$ occurs before the last bead on runner $\pt-1$ (which is no earlier than position $\pt-1$), so in fact runner $\pt-2$ is the smallest runner of $\al\pt$. However, $\la$ has a space on runner $\pt-1$ no later than position $-2p+\pt-1$, which means that $\al\pt$ has a space on runner $\pt+1$ no later than position $-2p+\pt+1$, i.e.\ earlier than the last bead on runner $\pt-2$. So by \cref{qsstuff}(\ref{pt1}) $\al\pt$ is not \qs, a contradiction.
\end{pf}

Now we consider the case $i=1$.

\begin{propn}\label{1wt1}
Suppose $p\gs5$, and that $\al1\neq\la$ while $\al0=\la$. Suppose furthermore that runner $1$ is the smallest runner in the abacus display for \al1. Then either $\nrd1$ or $\nrd{p-1}$.
\end{propn}

\begin{pf}
Since runner $1$ is the smallest in the abacus display for \al1, the last bead on this runner occurs before the space in position $\pt-p$, i.e.\ no later than position $1-p$. Hence for every $l\gs0$ at most one of the positions $lp$ and $lp+1$ is occupied in the abacus display for $\la$. If none of the positions $lp$ for $l\gs0$ is occupied, then we can just copy the proof of \cref{iwti} to deduce that \nrd1. So assume that position $kp$ is occupied for some $k\gs0$. In fact we must have $k>0$, since if position $0$ were occupied then $\la$ would have a removable $0$-node. By the above analysis position $kp+1$ is vacant, and position $1$ is occupied, since otherwise we would have $\quo0{(\al1)}\neq\varnothing$.

Since $\la$ is self-conjugate, this means that positions $-1$ and $-kp-2$ are occupied, while positions $-2$ and $-kp-1$ are vacant. Furthermore, for each $0>l>-k$ at least one of the positions $lp-2$ and $lp-1$ is occupied, and in addition at least one of the positions $-p$ and $1-p$ is occupied (again since $\quo0{(\al1)}=\varnothing$). So the number of beads between positions $-kp-1$ and $-2$ is at least $k$, and so \nrd{p-1}.
\end{pf}

\begin{rmk}
It turns out that the results in this subsection are to some extent redundant, because (given our standing assumptions) if $1<i<\pt-1$ and $\al j=\la$ for $0\ls j<i$, then $\al i=\la$ also. However, proving this statement appears to be just as difficult as the results we have proved in this section, which in any case are still needed for the cases $i=1,\pt-1$.
\end{rmk}

\subsection{Third case}\label{ptcase1sec}

We are left with the situation where $\al i=\la$ for $0\ls i<\pt$. We now have two possibilities to consider, since $\al\pt$ could be an \rtwo{} or a \jm. In this section we address the first of these possibilities.

\begin{propn}\label{ar}
Suppose that $\al j=\la$ for all $0\ls j<\pt$, and that \al\pt{} is an \rtwo. Then either \nrd\pt{} or \nrd{\pt+1}, or $p=3$ and $\la=(4^3,3)$.
\end{propn}

\begin{pf}
Because $\al\pt$ is an \rtwo, we have $\quo\pt{(\al\pt)}=(1)$, which means that in the abacus display for \al\pt{} (and hence in the abacus display for $\la$):
\begin{itemize}
\item
for $l=0$ and for every $l\ls-2$ there are at least two beads among positions $lp+\pt-1$, $lp+\pt$ and $lp+\pt+1$; and
\item
for $l=-1$ and for every $l\gs1$ there is at most one bead among positions $lp+\pt-1$, $lp+\pt$ and $lp+\pt+1$.
\end{itemize}
Note also that if there is a bead in position $\pt-1$ in the abacus display for $\la$, then since $\al j=\la$ for all $j<\pt$, there are beads in positions $\pt-2,\pt-3,\dots,-\pt$. But this contradicts the self-conjugacy of $\la$, so position $\pt-1$ is vacant, and hence positions $\pt$ and $\pt+1$ are occupied. Symmetrically, position $-\pt$ is occupied and positions $-\pt-2$ and $-\pt-1$ are vacant.

$\la$ must have at least one addable $\pt$-node, so there is an integer $k$ such that position $kp+\pt-1$ is occupied and position $kp+\pt$ is vacant. We consider separately the cases $k\ls-2$ and $k\gs1$.

If $k\ls-2$, then there are occupied positions $\pt$ and $kp+\pt-1$ and vacant positions $\pt-1$ and $kp+\pt$; from our observations so far, the number of intervening beads is at least $-2(k+1)\gs-k$, and so \nrd\pt.

Alternatively, suppose $k\gs1$. The condition that $\al j=\la$ for $0\ls j<\pt$ means that positions $kp+\pt-2,kp+\pt-3,\dots,(k-1)p+\pt+2$ are all occupied; by self-conjugacy, positions $-kp+\pt-2,-kp+\pt-3,\dots,(-k-1)p+\pt+2$ are all vacant. But now if $p\gs5$ then the first space on runner $0$ is earlier than the last bead on runner $p-1$ and vice versa, so that \al\pt{} is not \qs; contradiction.

We are left with the case where $p=3$ and $k\gs1$. In this case positions $-1$ and $-3k-2$ are occupied, while positions $-2$ and $-3k-1$ are vacant. The number of intervening beads is at least $2(k-1)$, so if $k\gs2$ then we have \nrd2. So we can assume that the only bead on runner $0$ with a space immediately after it is in position $3$. Now consider positions $6$, $7$ and $8$ on the abacus; from above, there is at most one bead among these three positions. There cannot be a bead in position $6$; in addition, there cannot be a bead in position $7$ (since then there would be a bead in position $-9$ and a space in position $-8$, again contrary to hypothesis). If there is a bead in position $8$, then we have beads in positions $8$ and $-5$ and spaces in positions $7$ and $-4$, with intervening beads in positions $3$, $2$, $1$ and $-1$, so \nrd2.

So we can assume that positions $6$, $7$ and $8$ are all vacant. If any position after position $8$ is occupied, then we have $\quo0{(\al1)}\neq\varnothing$; but since $p=3$, our assumptions on $\la$ means that \al1 has $3$-weight $1$, so $\quo0{(\al1)}=\varnothing$. So all positions after position $8$ are also vacant. The only possibilities left are that $\la=(4^3,3)$ or $(6,5^3,4,1)$, but in the latter case we again have \nrd2.
\end{pf}

\subsection{Fourth case}\label{ptcase2sec}

\newcommand\smab[4]{\raisebox{.5pt}{$\begin{array}{|@{\,}c@{\,}|@{\,}c@{\,}|}\hline&\\[-12pt]#1&\raisebox{-1pt}{$\sabacus(.85,#2#3#4)$}\\\hline\end{array}$}}

Now we come to the most difficult case. In addition to our standing assumptions, \emph{we assume throughout this subsection that $\al\pt$ is a \jm, and $\al j=\la$ for all $0\ls j<\pt$}. In particular, this means that $\quo\pt{(\al\pt)}=\varnothing$, which immediately gives some information about the middle three runners of the abacus display for $\la$.

\begin{lemma}\label{abstructure}
Suppose $k\in\bbz$. Then in the abacus display for $\la$:
\begin{itemize}
\item
at least two of the positions $kp+\pt-1$, $kp+\pt$ and $kp+\pt+1$ are occupied if $k<0$, while
\item
at most one of the positions $kp+\pt-1$, $kp+\pt$ and $kp+\pt+1$ are occupied if $k\gs0$.
\end{itemize}
\end{lemma}

In the next few results we analyse various possibilities for the configuration of the middle three runners of the abacus display for $\la$. To help us, we introduce some notation in which we give portions of the abacus diagram focusing on the middle three positions in a given row. For example, given an integer $k$ we may write \smab kbnb to mean that positions $kp+\pt-1$ and $kp+\pt+1$ are occupied in the abacus display for $\la$, while position $kp+\pt$ is vacant.

\begin{lemma}\label{7.1}
Suppose that \smab kbnn for some $k$. Then \nrd{-\pt} and \nrd{\pt(-\pt)}.
\end{lemma}

\begin{pf}
The fact that $\al j=\la$ for all $j<\pt$ means that positions $kp+\pt-2,kp+\pt-3,\dots,(k-1)p+\pt+1$ are occupied. In particular, this means that $k$ cannot equal $0$, since then by self-conjugacy we would have \smab{-1}bbn, a contradiction. So $k\gs1$, and hence \smab{k-1}nnb (since at most one of the positions $kp+\pt-1$, $kp+\pt$ and $kp+\pt+1$ is occupied).

The \scy{} of $\la$ now implies that \smab{-k-1}bbn and \smab{-k}nbb. In particular, there are beads in positions $(k-1)p+\pt+1$ and $-(k+1)p+\pt$ and spaces in positions $(k-1)p+\pt$ and $-(k+1)p+\pt+1$. Since $\la$ is \sc, exactly half of the positions from $-kp+\pt$ to $(k-1)p+\pt$ inclusive are occupied, so there are $\frac12((2k-1)p+1)$ beads in these positions; in particular, there are at least $2k$ beads between positions $-(k+1)p+\pt+1$ and $(k-1)p+\pt$, and so \nrd{-\pt}.

Exactly the same analysis applies to the partition $\reall{\pt}\la$; the only essential difference is that the bead in position $k-p+\pt$ has moved to position $-kp+\pt-1$, but this makes no difference to the count of intervening beads, so $\noregdown{-\pt}{\reall\pt\la}$, which gives \nrd{\pt(-\pt)}.
\end{pf}

\begin{lemma}\label{7.2.1}
Suppose there is no $k$ for which \smab kbnn. Then \smab knbn for some $k$, and hence \nrd\pt.
\end{lemma}

\begin{pf}
If there is no $k$ for which \smab knbn, then for every $k\gs0$ we have either \smab knnb or \smab knnn, and so for every $k<0$ we have either \smab knbb or \smab kbbb. But now $\la$ has no addable $\pt$-nodes; contradiction.

So \smab knbn for some $k$. So (by \scy) there are beads in positions $kp+\pt$ and $-(k+1)p+\pt-1$, and spaces in positions $kp+\pt-1$ and $-(k+1)p+\pt$. By a similar argument to that used in the proof of \cref{7.1}, the number of intervening beads is at least $2k+1$, and so \nrd\pt.
\end{pf}

\begin{lemma}\label{7.2.2}
If \smab knbn and \smab lnnb for some $k>l$, then \nrd{(\pt+1)\pt(\pt+1)}.
\end{lemma}

\begin{pf}
The self-conjugacy of $\la$ implies that \smab{-k-1}bnb and \smab{-l-1}nbb. So in the abacus display for $\reall{(\pt+1)\pt}\la$, positions $-(l+1)p+\pt+1$ and $-(k+1)p+\pt$ are occupied, while positions $-(l+1)p+\pt$ and $-(k+1)p+\pt+1$ are vacant. The number of intervening beads is at least $2(k-l-1)+1\gs k-l$, and so $\noregdown{-\pt}{\reall{(-\pt)\pt}\la}$, i.e.\ \nrd{(-\pt)\pt(-\pt)}.
\end{pf}

\begin{lemma}\label{7.2.3}
Suppose that for some $k<l$ we have \smab knbn and \smab lnnb, and that positions $(k-1)p+\pt+1$ and $(l+1)p+\pt$ are vacant, and define $\mu$ to be the partition obtained from $\la$ by repeatedly adding all addable nodes of residue not equal to $\pt$ or $-\pt$. Then $\noregdown\pt\mu$ and $\noregdown{\pt+1}\mu$, and hence $\spe\la$ has at least three \cf s.
\end{lemma}

\begin{pf}
Since position $(k-1)p+\pt+1$ is vacant and $\al j=\la$ for all $j<\pt$, positions $(k-1)p+\pt+2,\dots,kp+\pt-2$ are also vacant. Hence in the abacus display for $\mu$, position $kp+\pt-1$ is vacant. So the abacus display for $\mu$ has beads in positions $kp+\pt$ and $-(k+1)p+\pt-1$, and spaces in positions $kp+\pt-1$ and $-(k+1)p+\pt$. The number of intervening beads is $\frac12((2k+1)p-1)$, which is at least $2k+1$. So $\noregdown{\pt}\mu$.

Now we show that $\noregdown{-\pt}\mu$. Since position $(l+1)p+\pt$ is vacant in $\la$, it is also vacant in $\mu$. Since position $lp+\pt+1$ is occupied in $\la$, position $(l+1)p+\pt-1$ is occupied in $\mu$. As we have already seen, position $kp+\pt-1$ is vacant in $\mu$, while position $kp+\pt$ is occupied. $\mu$ is \sc, so the abacus display for $\mu$ has beads in positions $-(k+1)+\pt+1$ and $-(l+2)p+\pt$, and spaces in positions $-(k+1)+\pt$ and $-(l+2)p+\pt+1$. The number of intervening beads is at least $l-k+1$, and so $\noregdown{-\pt}\mu$.

So we have $\noregdown{\pt(-\pt)\pt}\mu$ and $\noregdown{(-\pt)\pt(-\pt)}\mu$, so by \cref{ligh0,ligh1} $\spe\mu$ has at least three \cf s. So by \cref{213f} $\spe\la$ does too.
\end{pf}

\begin{lemma}\label{7.3}
Suppose that \smab0nbn. Then there is some $l\gs1$ such that either \smab lbnn or \smab lnbn.
\end{lemma}

\begin{pf}
Suppose not, i.e.\ for every $l\gs1$ either \smab lnnb or \smab lnnn. We will show that $\la$ is \qs, contrary to assumption. Since $\al0=\la$, position $0$ is vacant in the abacus display for $\la$. Since $\al j=\la$ for all $j<\pt$, this then means that positions $1,\dots,\pt-2$ are also vacant. In addition, the fact that position $\pt+1$ is vacant means that positions $\pt+2,\dots,p+\pt-2$ are also vacant. Symmetrically, positions $-2p+\pt+2,\dots,-p+\pt-2$ and $-p+\pt+2,\dots,-1$ are all occupied.

Now we claim that there is at least one value of $l$ for which \smab lnnb. If not, then position $lp+\pt+1$ is vacant in $\la$ for all $l\gs0$, and hence (since $\al j=\la$ for $j<\pt$) so are positions $lp+\pt+2,\dots,(l+1)p+\pt-2$. In other words, all positions $k\gs0$ apart from position $\pt$ are vacant, and all positions $k<0$ apart from position $-p+\pt$ are occupied. Hence $\la$ is the partition $(\pt+1,1^\pt)$; but this has $p$-weight $1$, contrary to assumption.

By assumption $\al\pt$ is a \jm, and in particular is \qs. We claim that runner $\pt-1$ is the largest runner in \al\pt. Choose $l$ such that \smab lnnb. Then the abacus display for \al\pt{} has a bead in position $lp+\pt-1\gs p+\pt-1$, and every runner apart from runner $\pt-1$ has a space before this position. So by \cref{qsstuff}(\ref{pt1}) runner $\pt-1$ is largest. Symmetrically, runner $\pt+1$ is smallest, and hence for every $i\neq\pt-1,\pt,\pt+1$ we have $\quo i\la=\quo i{(\al\pt)}=\varnothing$. So the abacus display for $\la$ has the following properties.
\begin{itemize}
\item
For any $i\neq\pt-1,\pt,\pt+1$, the last bead on runner $i$ occurs in position $-p+i$, and the first space in position $i$.
\item
On runner $\pt-1$, the last bead occurs in position $-p+\pt-1$.
\item
On runner $\pt$, the last bead occurs in position $\pt$ and the first space in position $-p+\pt$.
\item
On runner $\pt+1$, the first space occurs in position $\pt+1$.
\end{itemize}
It follows that $\la$ is \qs; contradiction.
\end{pf}

\begin{lemma}\label{7.4}
Suppose that for every $k\gs0$ either \smab knbn or \smab knnn. Then there is a $p$-restricted partition $\nu\notin\{\la\rest,\mull p{\la\rest}\}$ such that $\ho{\spe\la}{\spe\nu}\neq0$. Hence $\spe\la$ has at least three \cf s.
\end{lemma}

\begin{pf}
Since position $lp+\pt+1$ is vacant in the abacus display for $\la$ for every $l\gs0$ and $\al i=\la$ for all $0\ls i<\pt$, positions $lp+\pt+2,\dots,(l+1)p+\pt-2$ are also vacant. In addition, position $0$ is vacant, and so positions $1,\dots,\pt-2$ are vacant as well. So every position $k\gs0$ not on runner $\pt$ is vacant, and symmetrically every position $k<0$ not on runner $\pt$ is occupied.

Now we observe that there must be at least two values of $l$ for which \smab lnbn: if there are no such values, then $\la=\varnothing$, while if there is only one such value then $\la$ is an \rone, and either way our assumptions are violated. Let $k$ be minimal such that \smab knbn, and define a partition $\mu$ by moving the bead in position $kp+\pt$ to position $kp+\pt-1$, and the bead in position $-(k+1)p+\pt-1$ to position $-(k+1)p+\pt$. Then $\mu$ is obtained from $\la$ by removing a node and adding a node lower down of the same residue. The number of beads in positions $-(k+1)p+\pt+1,\dots,kp+\pt-2$ is $kp+\pt\gs2k+1$, and so (arguing as in the proof of \cref{arab}) the added node is in a later ramp than the removed node. Furthermore, the minimality of $k$ means that there are no addable nodes (of any residue) between the removed node and the added node. Hence by \cref{mainhom} there is a non-zero homomorphism $\spe\la\to\spe\nu$, where $\nu=\mu\rest$. Since $\soc{\spe\nu}\cong\sid\nu$, this means that $\sid\nu$ occurs as a composition factor of $\spe\la$, and it just remains to show that $\nu\notin\{\la\rest,\mull p{\la\rest}\}$.

Since $\mu$ is obtained from $\la$ by moving a node to a different ramp, $\la$ and $\mu$ cannot possibly have the same restrictisation, so $\nu\neq\la\rest$. To show that $\nu\neq\mull p{\la\rest}$, we claim that $\len{(\la\rest)}=\len\la$. To see this, we examine the Young diagram of $\la$. Let $l_1>\dots>l_r$ be the values of $l$ for which \smab lnbn. Then $\la$ is the partition
\[
\big(l_1p+\pt+1,l_2p+\pt+2,\dots,l_rp+\pt+r,r^{l_rp+\pt},(r-1)^{(l_{r-1}-l_r)p-1},\dots,1^{(l_1-l_2)p-1}\big).
\]
The last non-empty ramp of $\la$ must contain a removable node $(a,b)$; we cannot have $a\ls r$, since then $a<b$ and the removable node $(b,a)$ lies in a later ramp; so $(a,b)$ is a node of the form $(l_ip+\pt+i,i)$, for $1\ls i\ls r$. This node lies in ramp $(l_i(p-1)+i)p+(p-1)\pt-p$, and this is obviously maximised for $i=1$. So the node in the last non-empty row of $\la$ lies in the last non-empty ramp, and hence $\len{{\la\rest}}=\len\la$.

A similar argument applies to $\mu$, and so we have $\len\nu=\len\mu=\len\la$, the latter equality following from the fact that there are least two values of $l$ for which \smab lnbn. Now observe that $\len\la\equiv\pt+1\ppmod p$, and so if $\nu=\mull p{\la\rest}$ then $\len{{\la\rest}}+\len{\mull p{\la\rest}}\equiv1\ppmod p$, contradicting \cref{mulllemma}.
\end{pf}

Now we combine the results of this subsection.

\begin{propn}\label{7.5}
Suppose $\al j=\la$ for all $j<\pt$, and that \al\pt{} is a \jm. Then $\spe\la$ has at least three \cf s.
\end{propn}

\begin{pf}
We consider the various possibilities for the middle three runners of the abacus display for $\la$. If we have \smab kbnn for some $k$, then by \cref{7.1} we have $\nrd{\pt(-\pt)\pt}$ and $\nrd{(-\pt)\pt(-\pt)}$, so we are done by \cref{ligh1}(\ref{lige}); so assume there is no such $k$. Now by \cref{7.2.1} we have \smab knbn for some $k$, and by \cref{7.3} we can assume $k>0$. If there is no $l$ such that \smab lnnb then we are done by \cref{7.4}, so we assume there is at least one such $l$. If $l<k$, then we are done by \cref{7.2.1,7.2.2} (using \cref{ligh1}(\ref{lige})), so we may assume that for every $k,l$ with \smab knbn and \smab lnnb we have $k<l$. Now taking $k,l$ maximal such that \smab knbn and \smab lnnb, positions $(k-1)p+\pt+1$ and $(l+1)p+\pt$ in the abacus display for $\la$ are vacant, and so we are done by \cref{7.2.3}.
\end{pf}

\subsection{The proof of \cref{main}}

Now we combine the results proved in this section to give a proof of \cref{main}. As noted in \cref{altsec}, the `if' part has already been proved in \cite{altred}, so we need only prove the `only if' part. In other words, we must prove that if $\la$ is a \sc{} partition not in \twofac{} or \jmset, then $\spe\la$ has at least three composition factors. If $\la$ is \qs, then \cref{mainqs} gives the result, so we can assume this is not the case. We proceed by induction on $|\la|$, using the partitions $\al i$. If there is any $i$ for which $\al i\neq\la$ and $\al i\notin\twofac\cup\jmset$, then by \cref{irredspecht} and by induction $\spe{\al i}$ has at least three composition factors, so by \cref{ligh1}(\ref{liga}) $\spe\la$ does too. In addition, if $i<\pt$ and $\spe{\al i}$ is irreducible then by \cref{ligh1}(\ref{ligb}) $\spe\la$ has at least three composition factors. So we can assume that $\al i\in\twofac\cup\{\la\}$ for each $0\ls i<\pt$, while $\al\pt\in\twofac\cup\jmset\cup\{\la\}$. So the assumptions listed in \cref{assumpsec} apply. $\la$ must have at least one removable node, so for some $i$ we have $\al i\neq i$. Let $0\ls i<\pt$ be minimal with this property, and consider the possibilities for $i$.

If $i=0$, then by \cref{0nowt,0somewt} we have either $\nrd0$, or $p\gs5$ and $\nrd{p-1}$. So we are done by \cref{ligh1}(\ref{ligc}).

If $1\ls i<\pt$, then by \cref{biggestorsmallest,iwti-1,1wt0,iwti,1wt1} we have either $\nrd i$ or $\nrd{p-i}$, and so again we are done by \cref{ligh1}(\ref{ligc}).

If $i=\pt$, then from the discussion in \cref{assumpsec}, $\al\pt$ is either an \rtwo{} or a \jm. In the first case, \cref{ar} gives either $\nrd\pt$ or $\nrd{-\pt}$, or $p=3$ and $\la=(4^3,3)$. If $\nrd\pt$ or $\nrd{-\pt}$, then we are done by \cref{ligh1}(\ref{ligc}), while if $p=3$ and $\la=(4^3,3)$, then we have $\nrd{121}$ and $\nrd{212}$, so \cref{ligh1}(\ref{lige}) gives the result (or we can just use the readily-available decomposition numbers for $\sss{15}$). Finally, we have the case where $i=\pt$ and $\al\pt$ is a \jm. This is dealt with in \cref{7.5}.

\cref{main} now follows by induction.

\section{Irreducible representations which remain irreducible modulo every prime}\label{everyprimesec}

We conclude this paper with a corollary of our main theorem, in which we classify the irreducible representations of $\aaa n$ that remain irreducible modulo every prime. The result is unsurprising.

\begin{thm}\label{everyprime}
Suppose $\psi$ is an ordinary irreducible character of the alternating group $\aaa n$ and that $\psi$ remains irreducible modulo every prime. Then $\psi$ is one-dimensional.
\end{thm}

The same result for the symmetric groups was proved by Kleshchev and Premet in \cite{kp}, though an easier proof \cite{jmp2} follows from the classification of irreducible Specht modules in characteristic $2$.

\begin{pf}[Proof of \cref{everyprime}]
Assume $n\gs2$, with the case $n=1$ being trivial, and use the notation of \cref{altsec}.

Suppose $\psi=\psi^\la$, for $\la\neq\la'$. Then $\psi$ is the restriction to $\aaa n$ of an irreducible character $\chi^\la$ of $\sss n$. By \cite[Proposition 2.11]{altred} $\psi^\la$ remains irreducible modulo a given prime $p$ if and only if $\chi^\la$ does, so the result follows from the corresponding result for the symmetric groups.

Now suppose $\psi=\psi^{\la\pm}$ for $\la=\la'$. By \cite[Theorem 3.1]{altred}, $\psi$ is reducible modulo $2$ unless $\la=(2^2)$ or $\la=(r,r-1,\dots,1)$ for some $r\gs2$. So by \cref{main}, $\psi$ remains irreducible modulo $2$ and modulo $3$ if and only if $\la=(2^2)$ or $(2,1)$. But in both of these cases $\psi$ is one-dimensional.
\end{pf}

\section{Index of notation}

For the reader's convenience we conclude with an index of the notation we use in this paper. We provide references to the relevant subsections.

\newlength\colwi
\newlength\colwii
\newlength\colwiii
\setlength\colwi{2.5cm}
\setlength\colwiii{1cm}
\setlength\colwii\textwidth
\addtolength\colwii{-\colwi}
\addtolength\colwii{-\colwiii}
\addtolength\colwii{-1em}

\subsubsection*{Basic objects}
\vspace{-\topsep}
\begin{longtable}{@{}p{\colwi}p{\colwii}p{\colwiii}@{}}
$\sss n$&the symmetric group on $\{1,\dots,n\}$&\\
$\aaa n$&the alternating group on $\{1,\dots,n\}$\\
$\bbf$&a field\\
$p$&the characteristic of $\bbf$ (taking $p=\infty$ if $\bbq\subseteq\bbf$)\\
$h$&$\frac12(p-1)$ (when $p$ is odd)\\
$\bbf\sss n\mood$&the category of $\bbf\sss n$-modules\\
$\res_HM$&the restriction of a module $M$ to a subgroup $H$\\
\end{longtable}

\subsubsection*{Partitions}
\vspace{-\topsep}
\begin{longtable}{@{}p{\colwi}p{\colwii}p{\colwiii}@{}}
$\la$&a composition or its Young diagram&\ref{partnsubsec}\\
$|\la|$&the size (i.e.\ the number of nodes) of a composition $\la$&\ref{partnsubsec}\\
$\varnothing$&the unique partition of $0$&\ref{partnsubsec}\\
$\la'$&the conjugate to a partition $\la$&\ref{partnsubsec}\\
$\la\rest$&the \emph{$p$-restrictisation} of a partition $\la$&\ref{partnsubsec}\\
$\lad l\la$&the number of nodes of $\la$ in ramp $l$&\ref{partnsubsec}\\
$\alad l\la$&the number of addable nodes of $\la$ in ramp $l$&\ref{partnsubsec}\\
$\rlad l\la$&the number of removable nodes of $\la$ in ramp $l$&\ref{partnsubsec}\\
$\mull p\la$& the image of a $p$-restricted partition $\la$ under the Mullineux map&\ref{signsubsec}\\
$\jmset$&the set of \jm s&\ref{jmsec}\\
$\twofac$&the set of self-conjugate partitions (conjecturally) labelling Specht modules of composition length $2$&\ref{altsec}\\
$\lir\gamma\alpha\beta$&the Littlewood--Richardson coefficient corresponding to partitions $\alpha,\beta,\gamma$&\ref{rouqsec}\\
\end{longtable}

\subsubsection*{Representations of the symmetric group}
\vspace{-\topsep}
\begin{longtable}{@{}p{\colwi}p{\colwii}p{\colwiii}@{}}
$\yper\la$&the Young permutation module indexed by a composition $\la$&\ref{partnsubsec}\\
$\spe\la$&the Specht module indexed by a partition $\la$&\ref{partnsubsec}\\
$\sid\la$&the simple module indexed by a $p$-restricted partition $\la$&\ref{partnsubsec}\\
$\sgn$&the one-dimensional sign representation of $\bbf\sss n$&\ref{signsubsec}\\
$\chi^\la$&the character of $\spe\la$, when $p=0$&\ref{altsec}\\
$\psi^\la$&the restriction of $\chi^\la$ to $\aaa n$, when $\la\neq\la'$&\ref{altsec}\\
$\psi^{\la+},\psi^{\la-}$&the irreducible summands of the restriction of $\chi^\la$ to $\aaa n$, when $\la=\la'$&\ref{altsec}
\end{longtable}

\subsubsection*{Restriction functors}
\vspace{-\topsep}
\begin{longtable}{@{}p{\colwi}p{\colwii}p{\colwiii}@{}}
$\e_i$&the $i$-restriction operator $\bbf\sss n\mood\to\bbf\sss{n-1}\mood$, for $i\in\{0,\dots,p-1\}$&\ref{eisubsec}\\
$\e_i^{(r)}$&the $r$th divided power of $\e_i$, for $r\gs0$&\ref{eisubsec}\\
$\eps iM$&$\max\rset{r\gs0}{\e_i^{(r)}M\neq0}$, for a non-zero module $M$&\ref{eisubsec}\\
$\e_i^{(\max)}M$&$\e_i^{(\eps iM)}M$&\ref{eisubsec}\\
$\remo i\la$&the number of removable $i$-nodes of a partition $\la$&\ref{eisubsec}\\
$\reall i\la$&the partition obtained from $\la$ by removing all the removable $i$-nodes&\ref{eisubsec}\\
$\reall{i_1,\dots,i_r}\la$&$\reall{i_r}{(\reall{i_2}{(\reall{i_1}\la)}\dots)}$&\ref{eisubsec}\\
$\al i$&the partition obtained from $\la$ by repeatedly removing all removable nodes of residue $\pm i$&\ref{rrsubsec}\\
$\addall i\la$&the partition obtained from $\la$ by adding all the addable $i$-nodes&\ref{eisubsec}\\
$\addall{\pm i}\la$&the partition obtained from $\la$ by repeatedly adding all addable nodes of residue $\pm i$&\ref{assumpsec}\\
$\nor i\la$&the number of normal $i$-nodes of a partition $\la$&\ref{eisubsec}\\
$\renor i\la$&the partition obtained from $\la$ by removing all the normal $i$-nodes&\ref{eisubsec}\\
$\noregdown{i_1,\dots,i_r}\la$&$\nor{i_l}{\renor{i_1\dots i_{l-1}}{(\la\rest)}}<\remo{i_l}{\reall{i_1\dots i_{l-1}}{\la}}$ for some $1\ls l\ls r$&\ref{rrsubsec}
\end{longtable}

\subsubsection*{The abacus}
\vspace{-\topsep}
\begin{longtable}{@{}p{\colwi}p{\colwii}p{\colwiii}@{}}
$\pqd\la$&the $p$-quotient of a partition $\la$&\ref{abdispsec}\\
$q_i(\la)$&the position of the first space on runner $i$ in the abacus display for the $p$-core of $\la$&\ref{rouqsec}\\
$\pi_\la$&the permutation of $\zpz$ such that $q_{\pi(0)}<\dots<q_{\pi(p-1)}$&\ref{rouqsec}\\
$\opqd\la$&the ordered $p$-quotient of $\la$&\ref{rouqsec}\\
$d_{\la\mu},a_{\la\mu}$&functions used in the formula for decomposable numbers for Rouquier partitions&\ref{rouqsec}\\
$\smab knbb$&positions $kp+\pt$ and $kp+\pt+1$ are occupied in the abacus display for $\la$, while position $kp+\pt-1$ is vacant&\ref{ptcase2sec}
\end{longtable}

\subsubsection*{Tableaux and homomorphisms}
\vspace{-\topsep}
\begin{longtable}{@{}p{\colwi}p{\colwii}p{\colwiii}@{}}
$T_{x,y}$&the entry in row $x$ and column $y$ of a tableau $T$&\ref{tabhomsubsec}\\
$\tht T$&the homomorphism $\yper\la\to\yper\mu$ labelled by a $\la$-tableau $T$ of type $\mu$&\ref{tabhomsubsec}\\
$\hth T$&the restriction of $\tht T$ to $\spe\la$&\ref{tabhomsubsec}\\
$\dom$&the dominance order on tableaux&\ref{tabhomsubsec}\\
$\la\stre$&a partition used in the construction of a magic $\la$-tableau&\ref{magicsubsec}\\
$\insm A$&a $\la$-tableau constructed from a $\la\stre$-tableau $A$&\ref{magicsubsec}\\
$\rtab\la$&the `restrictisation tableau' of shape $\la$ and type $\la\rest$&\ref{magicsubsec}
\end{longtable}

\subsubsection*{Multisets}
\vspace{-\topsep}
\begin{longtable}{@{}p{\colwi}p{\colwii}p{\colwiii}@{}}
$T^i$&the multiset of entries in the $i$th row of a tableau $T$&\ref{composubsec}\\
$X_i$&the multiplicity of a natural number $i$ as an element of a multiset $X$&\ref{composubsec}\\
$X\sqcup Y$&the multiset defined by $(X\sqcup Y)_i=X_i+Y_i$ for all $i$&\ref{composubsec}\\
$X+1$&the multiset obtained by adding $1$ to every element of $X$&\ref{composubsec}\\
$\{l\}^n$&the multiset with $n$ elements all equal to $l$&\ref{composubsec}
\end{longtable}

\end{document}